\documentclass[11pt,a4paper]{article}
\usepackage{amssymb}
\usepackage{amssymb}
\usepackage{amssymb}
\usepackage{amssymb}
\usepackage{amssymb}
\usepackage{amssymb}
\usepackage{amssymb}
\usepackage{amssymb}
\usepackage{amssymb}
\usepackage{amssymb}
\usepackage{amssymb}
\usepackage{amssymb}
\usepackage{graphicx}
\usepackage{amsmath}

\def\mineappendix{
        \setcounter{section}{1}
        \setcounter{subsection}{0}
        \def\thesection{\Alph{section}}
        \def\sectionap{\@startsection  {section}{1}{\z@}
                        {-3.5ex plus-1ex minus-.2ex} {0ex plus.2ex}
                        {\reset@font\Large\bf  Appendix:  \, }}}
\makeatother
\def\Proclaim #1. #2\par{\bigbreak\noindent{\sc#1.\enspace}{\it#2}\par}

\font\Bbbfont=msbm10
\newfam\msbfam
\textfont\msbfam=\Bbbfont \textfont\msbfam=\Bbbfont

\def\nn{\nonumber}

\newtheorem{Definition}{Definition}
\newtheorem{Theorem}{Theorem}%[section]
\newtheorem{Corollary}{Corollary}
\newtheorem{Proposition} {Proposition}
\newtheorem{Remark}{Remark}

\usepackage[dvipdfm, pdfstartview=FitH,
CJKbookmarks=true, bookmarksnumbered=true,
bookmarksopen=true,
citecolor=blue %¸Ä±äĿ¼¡¢ÒýÓõÄÑÕÉ«
,colorlinks=true
,pdfborder=1 %×¢Ê͵ô´ËÏîÔò½»²æÒýÓÃΪ²ÊÉ«±ß¿ò
,pdfstartpage=2 %´ò¿ªÊ±ÎªµÚ2Ò³
,pdfstartview=Fit]{hyperref} %´ò¿ªpdfʱΪ fit page
\allowdisplaybreaks
\title{\bf Generalized Bi-Schr\"odinger Flows and Vortex Filament on Symmetric Lie Algebras}
\author{Qing Ding$^{1}$\footnote{Email:qding@fudan.edu.cn; Fax: 86-21-65646073.}
\ \ \ and  \ \ \ Youde Wang$^{2}$\footnote{Email:wyd@math.ac.cn.}\\
\it $^1$ Institute of Mathematics and Key Lab. of Math. for Nonlinear Sciences  \\
\it Fudan University, Shanghai 200433, P.R. China\\
\it $^2$ Academy of Mathematics and Systematic Sciences\\
\it Chinese Academy Sciences, Beijing 100190, P.R. China}
\date{}
\begin{document}
\maketitle

\begin{abstract}
The theory of the vortex filament in three-dimensional fluid
dynamics, consisting mainly of the models up to the third-order
approximation (refer to \cite{AH,FM,FMo}), is an attractive subject
in both physics and mathematics. Many efforts have been devoted to
the extension of the theory to higher-dimensional symmetric Lie
algebras. However, such a generalization known in literature is
strongly restricted by the integrable method. In this article, we
endeavor to establish the third-order models of the vortex filament
on symmetric Lie algebras in a purely geometric way by generalized
bi-Schr\"odinger flows. Our generalization overcomes the limitation
of integrability and creates successfully the desired models on
Hermitian or para-Hermitian symmetric Lie algebras. Combining the
result in this article with what have been known in literature for
the leading-order and the second-order models, we actually exhibit
the basic models and the related theory of the vortex filament on
symmetric Lie algebras up to the third-order approximation.
\end{abstract}



\section{Introduction}
Vortex rings in fluid physics are usually regarded as invariant
states in the 3-dimensional inviscid incompressible fluid in which
the vortex lines are endowed with curvature that enables the rings
to propagate. The study of vortex rings can be traced back to
Helmholtz in 1858 (see \cite{Hel, Hel*}) and Kelvin in 1867 (see
\cite{Ke}). Since then, the theory of vortex filament dynamics in
inviscid fluid has become an attractive subject in both physics and
mathematics, and has been studied extensively and deeply by many
physicists and mathematicians (see \cite{AH, Da, Ha, FM, MGK}). We
should point out that, except for the hydrodynamic models in the
classical mechanics, one also needs to characterize the motion
behavior of vortex filament appearing in various physical models,
such as kinematics of interfaces in crystal growth \cite{Lan,BK},
viscous fingering in a Hele-Shaw cell \cite{ST}, charged fluid in a
neutralizing background \cite{UIY}, superconductors, superfluid
\cite{BJOS}, etc.

As an idealized model that essentially reveals some qualitative
behaviors of a curved vortex filament, the so called ``localized
induction equation (LIE)'' or in other words, ``localized induction
approximation (LIA)'' was first introduced in literature
\cite{AH}. The model in the Euclidean 3-space $\mathbb R^3$ is stated as follows.
The filament curve ${X}={X}(t,x)$, expressed by functions of the
arclength $x$ and the time $t$, evolves according to
\begin{eqnarray}
{X}_t=\kappa {\bf B},\label{Da}
\end{eqnarray}
where $\kappa$ is the curvature and ${\bf B}$ the binormal vector at
the position ${X}(t,x)$. This model is also called the binormal
motion of space curves in $\mathbb R^3$, derived by Luigi Sante
Da Rios, an Italian mathematician and
physicist, in 1906 (see \cite{Da}). The distinguishing features of LIE (\ref{Da}) are its
completely integrability and the fact that it is equivalent to the focus
nonlinear Schr\"odinger equation (NLS$^+$)
$$i\varphi_t+\varphi_{xx}+2|\varphi|^2\varphi=0$$
by Hasimoto transform (see \cite{Ha}). This produces the so-called
LIE-NLS correspondence (see \cite{LP} for example), that is, if a
curve $X$ evolves according to LIE (\ref{Da}), then the associated
complex function $\varphi$ given by $
\varphi(t,x)=\kappa(t,x)\exp\left(\int^x\tau(t,s)ds\right)$ evolves
according to NLS$^+$, where $\tau$ stands for the torsion curvature.
Magri \cite{Ma} unveiled the Hamiltonian structure behind the
integrability of LIE and created a recursion operator generating
successfully an infinite sequence of commuting vector fields. Based
on this, Langer and Perline in \cite{LP*} constructed the
counterparts for LIE. The resulting sequence of integrable vector
fields is called the ``localized induction hierarchy (LIH)'', the
first three terms of which are
\begin{eqnarray}
{\bf V}^{(1)}&=&\kappa {\bf B},~~~~{\bf
V}^{(2)}=\frac{1}{2}\kappa^2{\bf T}+\kappa_x{\bf N}+\kappa\tau {\bf
B},\nn\\
{\bf V}^{(3)}&=&\kappa^2\tau{\bf T}+(2\kappa_x\tau+\kappa\tau_x){\bf
N}+\left(\kappa\tau^2-\kappa_{xx}-\frac{1}{2}\kappa^3\right) {\bf
B},\nn\label{LIH}
\end{eqnarray} where $\{{\bf T},{\bf N},{\bf B}\}$ is the Frenet
frame along the filament curve. This further creates effective
mathematical methods and techniques to explore and describe the
property, hidden structure and symmetry of such vortex filaments.
Meanwhile, physically speaking, some other models are related to
particular adjustments of the velocity description of the vortex
filament curve with respect to ratio $\varepsilon$ of the core
radius to the ring radius (refer to \cite{FMo} and the references
therein for details). These models can be constructed by taking account
of the effect from the higher-order corrections to that of curvature
and torsion. For example, for the validity to $O(\varepsilon^2)$, a
vortex filament with axial velocity is reducible to a summation of
${\bf V}^{(1)}$ and ${\bf V}^{(2)}$, which was written down
explicitly by Fukumoto and Miyazaki in 1991 \cite{FM}. By omitting the
first term and re-normalizing, the Fukumoto-Miyazaki's model is
stated (see \S2 below for details) as
\begin{eqnarray}\frac{d X}{dt} =\frac{1}{2}\kappa^2{\bf T}+\kappa_x{\bf
N}+\kappa\tau{\bf B}={ X}_{xxx}+\frac{3}{2}{ X}_{xx}\times({
X}_x\times {X}_{xx}).\label{005}
\end{eqnarray}
Hence the Fukumoto-Miyazaki model is the second-order (physical)
correction model of the vortex filament in the Euclidean 3-space
${\mathbb R}^3$. According to this point of view, the LIE (\ref{Da})
is in fact the leading-order correction $O(\varepsilon)$-model of the
vortex filament in the Euclidean 3-space ${\mathbb R}^3$.

In 2000, the third-order correction $O(\varepsilon^3)$-model of the vortex
filament in the Euclidean 3-space ${\mathbb R}^3$ was derived by
Fukumoto and Moffatt (\cite{F,FMo}) (see \S2 below for more details)
and is written
\begin{eqnarray}
\frac{d X}{dt}=\lambda\bigg\{\kappa{\bf B}+\nu[\kappa^2\tau{\bf
T}+(2\kappa_x\tau+\kappa\tau_x){\bf
N}+(\kappa\tau^2-\kappa_{xx}){\bf B}]+\mu\kappa^3{\bf B}\bigg\}.
\label{I0}
\end{eqnarray}
However, the Fukumoto-Moffatt's model fails to be a summation of
${\bf V}^{(1)}$ and ${\bf V}^{(3)}$. Hence, surprisingly, the
Fukumoto-Moffatt's model of the vortex filament in the Euclidean 3-space
${\mathbb R}^3$ is in fact non-integrable and cannot be the fourth
member in the localized induction hierarchy in general (see
\cite{FMo}). All these models constitute a complete mathematical
treatment of the vortex filament up to the third-order approximation in
Riemannian 3-manifolds. The similar timelike and spacelike models of
vortex filament in the Minkowski 3-space $\mathbb R^{2,1}$ or
Lorentzian 3-manifolds, which possess totally different geometric
and dynamical properties by their timelike or spaclike stamp, are
constructed and characterized in a series of works in
\cite{D1,DI,dinghe,DWW,DWW1,DW1} (or see \S2 for details). These
models consist of interesting and indispensable supplements to the
models of the vortex filament in $\mathbb R^3$.

\medskip
On the other hand, in the past decades, many efforts have been
devoted to the extension of the vortex filament to higher
dimensional spaces. A successful generalization in literature seems
to be the ``localized induction (matrix) hierarchy'' on a symmetric
Lie algebra which depends crucially on the integrability (see, for
example, \cite{FK,LP,OS,TeUh}). Fordy and Kulish constructed in
\cite{FK} the matrix nonlinear Schr\"odinger equation associated to
a Hermitian symmetric space, which also promoted the integrable
matrix AKNS hierarchy on a Hermitian symmetric Lie algebra by the
spectral method. They showed that the matrix nonlinear Schr\"odinger
equation they constructed is in fact equivalent to a generalized
Heisenberg ferromagnet. Langer and Perline applied the technique of
Sym (\cite{Sym}) and Pohlmeyer \cite{Poh}) to produce geometric
realizations of the matrix nonlinear Schr\"odinger equation on
Hermitian symmetric Lie algebras (see \cite{LP}). In this process,
they illuminated the (integrable) ``localized induction (matrix)
hierarchy'' for arclength-parameterized curves evolving in a
Hermitian symmetric Lie algebra ${\bf g}$. The first three of them
are:
\begin{eqnarray}
{\widetilde\gamma}_t&=&{\widetilde\gamma}_x,\,\nonumber\\
{\widetilde\gamma}_t&=&[{\widetilde\gamma}_x,\,{\widetilde\gamma}_{xx}],\,\label{GLIE}\\
{\widetilde\gamma}_t&=&{\widetilde\gamma}_{xxx}+\frac{3}{2}[{\widetilde\gamma}_{xx},\,
[{\widetilde\gamma}_x,\,{\widetilde\gamma}_{xx}]],\label{GKDV}
\end{eqnarray}
where $[\cdot,\cdot]$ denotes the Lie bracket in ${\bf g}$.  We know
that, corresponding to the three types of models of the vortex
filament in $\mathbb R^3$ and $\mathbb R^{2,1}$, there are three
classes of important symmetric Lie algebras. The first class
consists of the Hermitian symmetric Lie algebras of compact type,
i.e. ${\bf g}=u(n)$ ($n\ge2$); the second class consists of the
Hermitian symmetric Lie algebras of noncompact type, i.e. ${\bf
g}=u(k,n-k)$ ($n\ge2$ and $1\le k\le n-1$); and the third class
consists of the para-Hermitian symmetric Lie algebras ${\bf
g}=gl(n,\mathbb R)$ ($n\ge2$) (refer to \S2 for more information).
When ${\bf g}=su(2),\, su(1,1)$ and $sl(2)$, Eq.(\ref{GLIE}) returns
respectively to the LIE (\ref{Da}) in $\mathbb R^3$, the timelike
and spacelike LIE in $\mathbb R^{2,1}$, and meanwhile, Eq.(\ref{GKDV})
goes respectively back to the Fukumoto-Miyazaki's model in $\mathbb
R^3$, and the corresponding timelike and spacelike models in
$\mathbb R^{2,1}$. We should point out that the Heisenberg
ferromagnet constructed by Fordy and Kulish in \cite{FK} is exactly
equivalent to model (\ref{GLIE}) of the localized induction
hierarchy on the Hermitian symmetric Lie algebra $u(n)$ ($n\ge2$).

\medskip
The geometric study of vortex filament in literature seems very
fascinating. It has been proved (refer to \cite{DW,D1,DI}, for
example) that the leading-order LIE model (\ref{Da}) in $\mathbb
R^3$, the timelike LIE (see (\ref{New041}) below) and the spacelike
LIE (see (\ref{New042}) below) in $\mathbb R^{2,1}$ of the vortex
filament are geometrically equivalent respectively to the
Schr\"odinger flow of maps from $\mathbb R$ to the 2-sphere $\mathbb
S^2\hookrightarrow {\mathbb R}^3$, the hyperbolic 2-space ${\mathbb
H}^2 \hookrightarrow {\mathbb R}^{2,1}$ and the de Sitter 2-space
${\mathbb S}^{1,1} \hookrightarrow {\mathbb R}^{2,1}$ and are also
respectively equivalent to the three typical 2nd-order integrable
systems of the AKNS hierarchy, i.e. the focus nonlinear
Schr\"odinger equation (NLS$^+$) (in this case, also see \cite{ZT}),
the defocus nonlinear Schr\"odinger equation (NLS$^-$)
$i\varphi_t+\varphi_{xx}-2|\varphi|^2\varphi=0$ and the nonlinear
heat equation (NLH) $q_t=q_{xx}-2qrq,~r_t=-r_{xx}+2rqr$ (for the
AKNS hierarchy, see \cite{AKNS,Ab0}). This indicates that the
LIE-NLS correspondence is valid for all the three NLS-type equations
and the LIE-type models. The second-order correction models of the
vortex filament in $\mathbb R^3$ and $\mathbb R^{2,1}$, i.e. the
Fukumoto-Miyazaki's model (\ref{005}) in $\mathbb R^3$ and its
sister timlelike and spacelike models (see \S2 for details) in
$\mathbb R^{2,1}$, are proved to be equivalent to the geometric KdV
flow of maps from $\mathbb R$ to $\mathbb S^2\hookrightarrow
{\mathbb R}^3$, ${\mathbb H}^2\hookrightarrow {\mathbb R}^{2,1}$ and
${\mathbb S}^{1,1}\hookrightarrow {\mathbb R}^{2,1}$ respectively
(see \cite{DWW,DW1}). This leads to that the four typical integrable
equations: $r_t+r_{xxx}+6rr_x=0$ (KdV), ~$r_t+r_{xxx}+6r^2r_x=0$
(mKdV), $\varphi_t=\varphi_{xxx} -{6}|\varphi|^2\varphi_x$
(cKdV$^-$) and $\varphi_t=\varphi_{xxx} +6|\varphi|^2\varphi_x=0$
(cKdV$^+$) existing in the 3rd-order integrable systems of the AKNS
hierarchy are interpreted in terms of the geometric KdV flows in a
unified way. Furthermore, it is also proved that the third-order
correction models of the vortex filament in $\mathbb R^3$ and
$\mathbb R^{2,1}$, i.e. the non-integrable Fukumoto-Moffatt's model
(\ref{I0}) in $\mathbb R^3$ and its brothers, the timelike and the
spacelike non-integrable models in $\mathbb R^{2,1}$, are equivalent
to the generalized bi-Schr\"odinger maps from $\mathbb R$ to
$\mathbb S^2\hookrightarrow {\mathbb R}^3$, ${\mathbb
H}^2\hookrightarrow {\mathbb R}^{2,1}$ and ${\mathbb
S}^{1,1}\hookrightarrow {\mathbb R}^{2,1}$ respectively (see
\cite{DWW1,DLW}). We must emphasize that, besides the concepts of
Schr\"odinger flows, the geometric KdV flows, the generalized
bi-Schr\"odinger flows, the K\"ahler structure on $\mathbb S^2$
(resp. $\mathbb H^2$) and the para-K\"ahler structure on $\mathbb
S^{1,1}$ play crucial roles in this aspect (see \cite{DI,DW1,DLW}
for details).

\medskip

Langer and Perline showed actually in \cite{LP} that the generalized
LIE-NLS correspondence is still valid on a Hermitian symmetric Lie
algebra ${\bf g}$. It was proved by Terng and Uhlenbeck in
\cite{TeUh} that the generalized LIE (\ref{GLIE}) on $u(n)$ is
equivalent to the Schr\"odinger flow of maps from $\mathbb R$ to the
Grassmannian manifold of compact type and by Chen in \cite{Chenb}
that the generalized LIE on $u(k,n-k)$ or $gl(n)$ are respectively
equivalent to the Schr\"odinger flow of maps from $\mathbb R$ to the
Grassmannian manifold of noncompact type or the para-Grassmannian
manifold, and are still equivalent respectively to the three typical
2nd-order matrix integrable systems of the matrix-AKNS hierarchy.
Along this route, Ding and He in \cite{dinghe} showed that the
third-order equation (\ref{GKDV}) of the localized induction
(matrix) hierarchy on the symmetric Lie algebra $sl(2n,\mathbb R)$
is equivalent to the geometric KdV flows from ${\mathbb R}$ to the
para-Grassmannian manifold, and hence the (integrable) matrix KdV
equation (see \cite{Lax,AF}) is created exactly by a motion of
Sym-Pohlmeyer curves in $sl(2n,\mathbb R)$ with initial data being
suitably restricted. Their result also signals that the third-order
equation (\ref{GKDV}) in the localized induction (matrix) hierarchy
on the symmetric Lie algebra $u(n)$ and $u(k, n-k)$ are respectively
the geometric KdV flows from ${\mathbb R}$ to the Grassmannian
manifold of compact type and of noncompact type. However, still
absent, but very important is to have the general third-order
non-integrable model of the vortex filament on symmetric Lie
algebras that corresponds to the (non-integrable) Fukumoto-Moffatt's
model in the Euclidean 3-space $\mathbb R^3$, and the timelike and
spacelike third-order correction models of the vortex filament in
the Minkowski 3-space $\mathbb R^{2,1}$.

\medskip
The purpose of this paper is just to establish the third-order
non-integrable models of the vortex filament on symmetric Lie
algebras such that they are natural extensions of the
Fukumoto-Moffatt's model (\ref{I0}) in $\mathbb{R}^3$.  Because of
the non-integrability, we can not employ integrable methods to
address the problem. In order to overcome the difficulty, we need to
carefully make correct characterizations of the third-order
Fukumoto-Moffatt's (physical) correction model of the vortex
filament motion in the Euclidean 3-space ${\mathbb R}^3$ in the way
they fit symmetric Lie algebras. The main idea applied here is the
K\"ahler or para-K\"ahler intrinsic characterization of the
symmetric spaces and then we derive the model on a symmetric Lie
algebra purely from the viewpoint of differential geometry. Recall
that we have introduced the concept of the so-called ``generalized
bi-Schr\"odinger flows" from a Riemannian manifold to K\"ahler and
para-K\"ahler manifolds in \cite{DWW1}. However, this concept is not
relevant to symmetric spaces in general. So our first task is to
correctly modify the geometric concept of generalized
bi-Schr\"odinger flows such that it is applicable to Hermitian and
para-Hermitian symmetric spaces. A key observation we obtained is
that the generalized bi-energy functional for maps $u:\mathbb R\to
(N,I,h)$ should be modified into
$$E_{\alpha,\beta, \gamma}(u)=\alpha\int_{\mathbb{R}}|\nabla u|^2dx
+\beta\int_{\mathbb{R}}|d\nabla u|^2dx+\gamma
\int_{\mathbb{R}}\langle R(\nabla u, I_u\nabla u)I_u\nabla u, \nabla
u\rangle dx,$$ where $\alpha$, $\beta$ and $\gamma$ are real
parameters, $(N, I, h)$ is a K\"ahler (or para-K\"ahler) manifold
with $I$ being the complex (or para-complex) structure,
$R(\cdot,\cdot)$ is the curvature operator of $N$ and
$\langle\cdot,\cdot\rangle$ stands for the inner product induced
from the metric $h$ on $N$. Then, the generalized bi-Schr\"odinger
flow in \cite{DWW1} is modified to be the Hamiltonian (or
para-Hamiltonian) gradient flow of the above generalized bi-energy
functional. Obviously, such a flow relies seriously on the intrinsic
geometry of the symmetric Lie algebras, as we shall elucidate in \S2
below.

By using this modification of the generalized bi-Schr\"odinger
flows, we successfully create the third-order non-integrable model
of the vortex filament on a symmetric Lie algebra ${\bf g}$ by the
equation of the generalized bi-Schr\"odinger flows from $\mathbb R$
to the corresponding symmetric space. Furthermore, by applying the
geometric concept of PDEs with the given curvature representation
which was proposed in the category of Yang-Mills theory in
\cite{ding1,ding2}, and with the aid of gauge transformations which
can be regarded as the generalization of the Hasimoto transform in
higher dimensions for non-integrable PDEs, we transform the
third-order model of the vortex filament on a symmetric Lie algebra
into a second to fourth order nonlinear Schr\"odinger-like
differential-integral matrix equation. All the results obtained in
this article coincide exactly with what we knew for the vortex
filament in $\mathbb R^3$ and $\mathbb R^{2,1}$ when the symmetric
Lie algebra ${\bf g}$ goes back to $su(2),\, su(1,1)$ and $sl(2)$.
Therefore, combining our result with what have been known in
literature for the leading and second order models, we actually
exhibit the basic models and the related theory of the vortex
filament on symmetric Lie algebras up to the third-order
approximation. The main results are summarized in Theorems 1-3.
\medskip

The article is organized as follows. In section 2, we give the
preliminaries about symmetric Lie algebras, the symmetric
Grassmannian or para-Grassmannian manifolds, the isospectral way to
obtain the matrix AKNS systems  and the geometric concept of
generalized bi-Schr\"{o}dinger maps from a Riemannian manifold to a
K\"ahler or para-K\"ahler manifold.  In section 3, we deduce the
equation of the generalized bi-Schr\"odinger maps from $\mathbb{R}$
to the three different types of the symmetric Grassmannian or
para-Grassmannian manifolds. In section 4,  the non-integrable
third-order models of the vortex filament on three different types
of symmetric Lie algebras are proved to be geometrical realizations
of Sym-Pohlmeyer curves in their own symmetric Lie algebras and are
transformed respectively to the second to fourth order nonlinear
Schr\"odinger-like differential-integral matrix equations by gauge
transformations. Conclusions and remarks are given in section 5.

\section{Preliminaries }
In this section, we will briefly mention the background of the vortex
filament in three dimensional Euclidean and Minkowski spaces. We also
recall some fundamental knowledge and facts about Hermitian or
para-Hermitian symmetric spaces and symmetric Lie algebras. Besides,
we need to elucidate the way of modifying the geometric
concept on generalized bi-Schr\"odinger flows.

\subsection{Three dimensional motions of vortex filaments}

Suppose that we have a fluid in a Riemannian 3-manifold $M^3$ or a
pseudo-Riemannian 3-manifold $M^{2,1}$ that evolves according to a
given one-parameter family of diffeomorphism yielding the position
of a fluid particle. Mathematically, the corresponding vortex
filament $\cal{L}$ is assumed to be parameterized by arclength $x$
and hence, we express a point (position) $\in \cal{L}$ as $X(t,x)$.
At a point $X(t,x)$ on the filament $\cal{L}$, we have the Frenet
frame $\{{\bf T},{\bf N},{\bf B}\}$, the curvature $k$ and the
torsion $\tau$. The kinematics is now described by
\begin{eqnarray}
v=v_T\varepsilon_1{\bf T}+v_1 \varepsilon_2{\bf
N}+v_2\varepsilon_3{\bf B}, \label{New02}\end{eqnarray}where
$\varepsilon_i$ ($i=1,2,3$) are respectively the sign of the tangent
vector ${\bf T}$, the normal vector ${\bf N}$ and the binormal
vector ${\bf B}$ which are called the causal characters of
$X_t(x)=X(t,x)$, $v=\frac{d X}{dt}(t,x)$ is the velocity of the
filament, and $v_T$, $v_1$ and $v_2$ are given quantities according
to the  physical background (these quantities cannot be arbitrarily
given, refer to \cite{DLW}). The general theory of the vortex
filament is to solve Eq.(\ref{New02}) (i.e. to determine $X(t,x)$)
for given $v_T$, $v_1$ and $v_2$.  It is obvious to see that
Eq.(\ref{New02}) is not easy to solve, so it is difficult to
understand the dynamical behavior of solutions to Eq.(\ref{New02}).

\medskip

The so-called ``localized induction equation'' (LIE), or in other
words, ``localized induction approximation'' (LIA), in $\mathbb R^3$
is an important idealized model of the vortex filament in
approximation put forward by Arms and Hama (see \cite{AH}), which
consists of solutions to the system (\ref{New02}) with
$\varepsilon_1=\varepsilon_2=\varepsilon_3=1$, assuming that
$v_T=v_1=0$ and $v_2=\kappa$ (refer to \cite{AH}). This model is
exactly the model (\ref{Da}), as has been mentioned previously in
the Introduction, which is the leading-order correction model of the
vortex filament. The LIE in Minkowski (or Lorentzian) 3-space
${\mathbb R}^{2,1}$ consists of solutions to the system
(\ref{New02}), assuming that $v_T=v_1=0$ and
$v_2=\varepsilon_2\kappa$. Hence, the corresponding timelike LIE
model in Minkowski 3-space $\mathbb R^{2,1}$ ($\varepsilon_1=-1$,
$\varepsilon_2=\varepsilon_3=1$) reads
\begin{eqnarray}
\frac{d X}{dt}~=\kappa {\bf B}. \label{New041}
\end{eqnarray}
One of the spacelike LIE models in Minkowski 3-space $\mathbb
R^{2,1}$ ($\varepsilon_1=1$, $\varepsilon_2=-1$ and
$\varepsilon_3=1$) states
\begin{eqnarray}
\frac{d X}{dt}~=-\kappa {\bf B}. \label{New042}
\end{eqnarray}
Another one is the case that $\varepsilon_1=\varepsilon_2=1$ and
$\varepsilon_3=-1$. One may verify directly that the system in this
case is equivalent to the system (\ref{New042}), so we only
consider in this article the case: $\varepsilon_1=1$,
$\varepsilon_2=-1$ and $\varepsilon_3=1$ when the filament is
spacelike. Eqs.(\ref{New041},\ref{New042}) are called respectively
the binormal motions of timelike and spacelike curves in $\mathbb
R^{2,1}$ in \cite{DI}.

\medskip

The second-order correction model of the vortex filament in $\mathbb R^3$
was deduced in 1991 by Fukumoto and Miyazaki (see \cite{FM}). This
model consists of solutions to (\ref{New02}) with
$\varepsilon_1=\varepsilon_2=\varepsilon_3=1$, assuming that
$v_T=\frac{\lambda}{2}\kappa^2$, $v_1=\lambda\kappa_{x}$ and
$v_2=\kappa+\lambda\kappa \tau$ ($\lambda$ is a real parameter),
i.e.,
\begin{eqnarray}
\frac{d X}{dt}&=&\kappa{\bf B}+\lambda\left(\frac{1}{2}\kappa^2{\bf
T}+\kappa_x{\bf N}+\kappa\tau{\bf B}\right), \nn\\
 &=&{
X}_x\times { X}_{xx}+\lambda\left[{ X}_{xxx}+\frac{3}{2}{
X}_{xx}\times({ X}_x\times { X}_{xx})\right].\label{New7}
\end{eqnarray}
Omitting the first term and re-normalizing, we obtain the
Fukumoto-Miyazaki's model (\ref{005}) indicated in the Introduction.
Similar to the Fukumoto-Miyazaki's model (\ref{005}), the
second-order timelike and spacelike correction models in the
Minkowski 3-space $\mathbb R^{2,1}$ are respectively deduced and
characterized in \cite{DW1} as follows:
\begin{eqnarray} \frac{d X}{dt}=\frac{1}{2}\kappa^2{\bf T}-\kappa_x{\bf
N}-\kappa\tau{\bf B}
 ={ X}_{xxx}+\frac{3}{2}{ X}_{xx}\dot\times({
X}_x\dot\times { X}_{xx}) \label{ab1}
\end{eqnarray}
and
\begin{eqnarray} \frac{d X}{dt}=\frac{1}{2}\kappa^2{\bf T}+\kappa_x{\bf
N}+\kappa\tau{\bf B}={ X}_{xxx}-\frac{3}{2}{X}_{xx}\dot\times({
X}_x\dot\times {X}_{xx}).\label{ab2}
\end{eqnarray}

\medskip

The third-order correction model of the vortex filament in the Euclidean
3-space ${\mathbb R}^3$ was derived by Fukumoto and Moffatt in 2000
(\cite{F,FMo}). According to the present description, the
Fukumoto-Moffatt's model is just Eq.(\ref{New02}) with
$\varepsilon_1=\varepsilon_2=\varepsilon_3=1$, assuming that
$v_T=\lambda\nu \kappa^2$, $v_1=\lambda\nu
(2\kappa_{x}\tau+\kappa\tau_{x})$ and
$v_2=\lambda\kappa+\lambda\nu(\kappa\tau^2-\kappa_{xx})+\lambda\mu
\kappa^3$ ($\lambda, \nu$ and $\mu$ are real parameters), which is
explicitly expressed by Eq.(\ref{I0}) stated in the Introduction.
Meanwhile, the third-order correction model in the Minkowski 3-space
$\mathbb R^{2,1}$ was proposed in \cite{DLW}, which consists of
solutions to (\ref{New02}) adjusted by assuming ($\lambda,
\nu$ and $\mu$ are also real parameters) that
$$v_T=\lambda \nu \varepsilon_3\kappa^2\tau,~v_1=\lambda\nu \varepsilon_2
\varepsilon_3(2\kappa_x\tau+\kappa\tau_x),~v_2= \lambda
\varepsilon_2\kappa-\lambda\nu\varepsilon_3
(\varepsilon_2\kappa_{xx}-\varepsilon_3\kappa\tau^2)+\lambda\mu\kappa^3.$$
Hence, the timelike third-order model of the vortex filament in $\mathbb
R^{2,1}$ ($\varepsilon_1=-1$ and $\varepsilon_2=\varepsilon_3=1$)
reads
\begin{eqnarray}
{X}_t=\lambda\bigg\{\kappa{\bf B}+\nu\left[-\kappa^2\tau{\bf
T}+(2\kappa_x\tau+\kappa\tau_x){\bf
N}-(\kappa_{xx}-\kappa\tau^2){\bf B}\right]+\mu\kappa^3{\bf
B}\bigg\} \label{I1}
\end{eqnarray}
and the spacelike one ($\varepsilon_1=1$, $\varepsilon_2=-1$ and
$\varepsilon_3=1$) is
\begin{eqnarray}
{X}_t=\lambda\bigg\{-\kappa{\bf B}+\nu\left[\kappa^2\tau{\bf
T}+(2\kappa_x\tau+\kappa\tau_x){\bf
N}+(\kappa_{xx}+\kappa\tau^2){\bf B}\right]+\mu\kappa^3{\bf
B}\bigg\}. \label{I2}
\end{eqnarray}

Models (\ref{Da}), (\ref{005}) and (\ref{I0}) are respectively the
leading order, the second order and the third order correction
equations of the vortex filament in the Euclidean 3-space $\mathbb
R^3$. These models characterize the phenomenon of the vortex
filament in a 3-dimensional Riemannian fluid manifold. The
exploitation of dynamical and geometric properties of these models
constitutes the complete theory of the vortex filament in $\mathbb
R^3$ up to the third-order approximation. The models (\ref{New041}),
(\ref{ab1}) and (\ref{I1}) for the timelike vortex filament  and the
models (\ref{New042}), (\ref{ab2}) and (\ref{I2}) for the spacelike
vortex filament in $\mathbb R^{2,1}$ provide the basic differential
equations of the theory of the vortex filament in a 3-dimensional
Lorentzian fluid manifold up to the third-order approximation.

\subsection{Symmetric Lie algebras and matrix-AKNS hierarchy}

A so-called symmetric Lie algebra ${\bf g}$ is a Lie algebra that
has a decomposition as a vector space sum: ${\bf g}={\bf
k}\oplus{\bf m}$ satisfying the (bracket) symmetric conditions:
$[{\bf k},{\bf k}]\subset {\bf k}$, $[{\bf m},{\bf m}]\subset {\bf
k}$ and $[{\bf k},{\bf m}]\subset {\bf m}$ (see \cite{AF,FK,He,LP}).
In such a symmetric Lie algebra there is an element denoted by
$\sigma_3$ in ${\bf k}$ such that ${\bf k}=\hbox{\rm
Kernel}(ad_{\sigma_3})=\{\chi\in{\bf g}|[\chi,\sigma_3]=0\}$. In
this section, let's first recall the three typical classes of
symmetric Lie algebras. Then we briefly review the matrix AKNS
hierarchy on a symmetric Lie algebra .

The first class of symmetric Lie algebras consists of Hermitian
symmetric Lie algebras $u(n)$ ($n\ge2$) with index $k$ ($1\le k< n$)
of compact type. In fact, for any given $1\le k<n$, let
\begin{eqnarray}
\sigma_3=\frac{\sqrt{-1}}{2}\left(\begin{array}{cc}
I_k&0\\0&-I_{n-k}\end{array}\right),\label{1}\end{eqnarray} and we see
that $u(n)$ is decomposable as $u(n)={\bf k}\oplus {\bf m}$
satisfy the symmetric conditions: $[{\bf k},{\bf k}]\subset {\bf
k}$, $[{\bf m},{\bf m}]\subset {\bf k}$, $[{\bf k},{\bf m}]\subset
{\bf m}$, where
$$
{\bf k}=\hbox{\rm
Kernel}(ad_{\sigma_3})=\left\{\left(\begin{array}{cc} A_{k\times
k}&0\\0&B_{(n-k)\times (n-k)}\end{array}\right)\in u(n)\right\}$$
and
$${\bf m}=\left\{\left(\begin{array}{cc} 0&U_{k\times
(n-k)}\\-U^*_{k\times (n-k)}&0\end{array}\right)\in u(n)\right\},
$$
where $U^*_{(n-k)\times k}$ stand for the transposed conjugate
matrix of $U_{k\times (n-k)}$. This indicates that $u(n)$ is a
symmetric Lie algebra depending on $k$ ($1\le k< n$). We should
point out that when $n=2m$ is even and $k=m$ the half of $n$, the
Hermitian symmetric Lie algebras $u(2m)$ with index $m$ may be
replaced exactly by $su(2m)$ with index $m$, since in this case,
$\sigma_3$ given by (\ref{1}) belongs to $su(2m)$ itself.

\medskip

The second class consists of Hermitian symmetric Lie algebras
$u(k,n-k)$ with index $k$ ($1\le k< n$) of noncompact type. In this
case, let $\sigma_3$ be also given by (\ref{1}) and we see that
$u(k,n-k)$ is decomposable as $u(k,n-k)={\bf k}\oplus {\bf m}$
satisfy the symmetric conditions, where
$$
{\bf k}=\hbox{\rm
Kernel}(ad_{\sigma_3})=\left\{\left(\begin{array}{cc} A_{k\times
k}&0\\0&B_{(n-k)\times (n-k)}\end{array}\right)\in
u(k,n-k)\right\}$$ and
$${\bf m}=\left\{\left(\begin{array}{cc} 0&U_{k\times
(n-k)}\\U^*_{(n-k)\times k}&0\end{array}\right)\in u(k,n-k)\right\}.
$$
Similarly, when $n=2m$ is even and $k=m$ the half of $n$, the
Hermitian symmetric Lie algebras $u(m,m)$ with index $m$ is replaced
just by $su(m,m)$ with index $m$. This is because $\sigma_3$ given
by (\ref{1}) belongs to $su(m,m)$ in this case.

\medskip

The third class consists of para-Hermitian symmetric Lie algebras
$gl(n,\mathbb R)$ ($n\ge2$) with index $k$ ($1\le k< n$). In this
case, for any given $1\le k<n$, setting
\begin{eqnarray}\sigma_3=\frac{{1}}{2}\left(\begin{array}{cc}
I_k&0\\0&-I_{n-k}\end{array}\right),\label{003}\end{eqnarray} we see
that $gl(n,{\mathbb R})={\bf k}\oplus {\bf m}$ is a symmetric Lie
algebra, where
$$
{\bf k}=\hbox{\rm
Kernel}(ad_{\sigma_3})=\left\{\left(\begin{array}{cc} A_{k\times
k}&0\\0&B_{(n-k)\times (n-k)}\end{array}\right)\in gl(n,\mathbb
R)\right\}$$ and
$${\bf m}=\left\{\left(\begin{array}{cc} 0&U^+_{k\times
(n-k)}\\U^-_{(n-k)\times k}&0\end{array}\right)\in sl(n,\mathbb
R)~\bigg|~\forall ~~U^+_{k\times (n-k)}~\hbox{and}~U^-_{(n-k)\times
k}\right\}.
$$
When $n=2m$ and $k=m$, the para-Hermitian symmetric Lie algebra
$gl(2m,\mathbb R)$ with index $m$ is replaced by the symmetric Lie
algebra $sl(2m,\mathbb R)$ with index $m$ (see \cite{dinghe}). It
should also be pointed out that para-Hermitian symmetric Lie
algebras relate to the geometry of the para-complex structure and the
para-K\"ahler structure, which has been studied since the first
papers of Rashevskii \cite{Ras} and Libermann \cite{Liber} until
now. For developments on this subject, referred to works
\cite{AMT,CFG} and the references therein.

We should mention that, for a given symmetric Lie algebra ${\bf
g}={\bf k}\oplus {\bf m}$ belonging to the above three classes and a
matrix $E\in{\bf g}$, we usually denote by $E^{(\hbox{diag})}$
(resp. $E^{(\hbox{off-diag})}$) its ${\bf k}$-part (resp. ${\bf
m}$-part) and is called the diagonal (resp. off-diagonal) part of
$E$.

\medskip

It is also well-known that the above three classes of symmetric Lie
algebras relate to three different types of symmetric spaces, which
may be presented by adjoint obits (refer to, for example,
\cite{BD,He}). The first one consists of the complex compact
(K\"ahler) Grassmannian manifolds:
\begin{eqnarray}
{G}_{n,k}=\left\{E^{-1}\sigma_3 E~|~ \forall ~E\in
U(n)\right\}:=U(n)/U(k)\times U(n-k), \label{GL1}
\end{eqnarray}
where $\sigma_3$ is given by (\ref{1}). The second one consists of
the complex noncompact (K\"ahler) Grassmannian manifolds:
\begin{eqnarray}
{G}_{n}^k=\left\{E^{-1}\sigma_3 E~|~ \forall ~E\in
U(k,n-k)\right\}:=U(k,n-k)/U(k)\times U(n-k), \label{GL2}
\end{eqnarray}
where $\sigma_3$ is still given by (\ref{1}). The third one consists
of the (para-K\"ahler) para-Grassmannian manifolds:
\begin{eqnarray}
{\widetilde G}_{n,k}=\left\{E^{-1}\sigma_3 E~|~ \forall ~E\in
GL(n,\mathbb R)\right\}:=GL(n,\mathbb R)/GL(k,\mathbb R)\times
GL(n-k,\mathbb R), \label{GL3}
\end{eqnarray}
where $\sigma_3$ is given by (\ref{003}). We point out that, as
mentioned above, we especially have
\begin{eqnarray}
{G}_{2m,m}&=&SU(2m)/SU(m)\times SU(m),\nonumber\\
{G}_{2m}^m&=&SU(m,m)/SU(m)\times SU(m)\nonumber
\end{eqnarray}
and
$$
{\widetilde G}_{2m,m}=SL(2m,\mathbb R)/SL(m,\mathbb R)\times
SL(m,\mathbb R).$$ For $M={G}_{n,k}$, ${G}_{n}^k$ and ${\widetilde
G}_{n,k}$, the K\"ahler structure $J_{u}$ (resp. the para-K\"ahler
structure $K_{u}$) at the point $u\in M$ is given by (see
\cite{TeUh,Chenb,dinghe})
\begin{eqnarray}
J_{u}~(\hbox{resp.}~K_{u})~=[u,\, \cdot]: ~T_{u}M\to
T_{u}M.\label{parastr}
\end{eqnarray}

\medskip
Next, in order to describe the matrix AKNS hierarchy, we consider an
overdetermined linear differential system as follows
\begin{eqnarray}
\varphi_x=(\lambda \sigma_3+Q)\varphi,~~~~\varphi_t=V\varphi.\label{3}
\end{eqnarray}
This system involves two independent variables $x$ (position) and
$t$ (time), and a spectral parameter $\lambda$. The eigenfunction
$\varphi$ takes values in a Lie group ${\bf G}$, while $(\lambda
\sigma_3+Q)$ and $V$ take values in the Lie algebra ${\bf g}$ of
${\bf G}$. For our purpose, ${\bf g}$ is required to be a symmetric
Lie algebra, and $\sigma_3$, $Q$ are suitably chosen for this
structure. In such a symmetric Lie algebra, $\sigma_3$ is chosen as
indicated above, $Q$ is required to be an ${\bf m}$-potential (i.e,
$Q(t,x)\in {\bf m}$ for all $t$) and $V$ a function of $Q$ and its
derivatives. One notes that the system (\ref{3}) gives the zero
curvature condition
\begin{eqnarray}
Q_t=V_x-[\lambda \sigma_3+Q,V].\label{4}\end{eqnarray} So, in
finding $V$ in terms of $Q$ and its derivatives such that the
integrability condition (\ref{4}) is satisfied, we invoke a
polynomial ansatz to $V$ as $V=\sum_{i=0}^NP^{(i)}\lambda^i$, where
$\forall i,$ $P^{(i)}$ is a function of $Q$ and its derivatives
independent of $\lambda$. The strategy here is to carry out $V$ for
a given integer $N$ and finally, to obtain a nonlinear integrable
PDE for the ${\bf m}$-potential $Q$. For details, refer to
\cite{FK,AF,LP}.

It is shown in \cite{FK} that the above process for $N=2$ on the
Hermitian symmetric Lie algebra $u(n)$, which is the so-called
second-order isospectral flow of the matrix Schr\"odinger equation
on $u(n)$, leads to the matrix nonlinear Schr\"odinger equation:
$$iQ_t+Q_{xx}+2QQ^*Q=0.$$
To obtain the matrix KdV equation, we work on the Lie algebra
${\bf g}=sl(2n,{\mathbb R})$. By a simple argument, the above
process for $N=3$ (i.e, the third-order isospectral flow) on
$sl(2n,\mathbb R)$ gives the ansatz expression for
$$
V=Q_{xx}-2Q^3+[Q,Q_x]+\left([\sigma_3,Q_x]+\frac{1}{2}[Q,[\sigma_3,Q]]\right)\lambda
+Q\lambda^2+\sigma_3\lambda^3,$$ which leads to the following
equation for an ${\bf m}$-potential $Q=\left(\begin{array}{cc}
0&U^+\\U^-&0\end{array}\right)$:
\begin{eqnarray}
Q_t=Q_{xxx}-2(Q^3)_x-[Q,[Q,Q_x]]. \label{401}
\end{eqnarray}
Eq.(\ref{401}) is rewritten as the following (noncommutative)
coupled matrix KdV equation:
\begin{eqnarray}
\left\{\begin{array}{c} U^+_t=U^+_{xxx}-3U^+U^-U^+_x-3U^+_xU^-U^+,\\
U^-_t=U^-_{xxx}-3U^-U^+U^-_x-3U^-_xU^+U^-.\end{array}\right.\label{5}
\end{eqnarray}
Generally, $U^+$ and $U^-$ in the coupled matrix KdV equation
(\ref{5}) are a pair of $k\times n$ matrix and $n\times k$ matrix.
Here we only deal with Eq.(\ref{5}) with $U^+$ and $U^-$ being
$n\times n$ matrices. Now by the reduction: $U=U^+$ and $U^-=I_n$,
we obtain the matrix KdV equation
\begin{eqnarray}
U_t=U_{xxx}-3UU_x-3U_xU \nn
\end{eqnarray}
from Eq.(\ref{5}). By a different reduction:
$U^-(t,x)=U^+(t,x)=U(t,x)$, Eq.(\ref{5}) leads to the modified
matrix KdV equation:
\begin{eqnarray}
U_t=U_{xxx}-3U^2U_x-3U_xU^2. \nn\label{501}
\end{eqnarray}

If one wants to obtain the fourth order integrable matrix
Schr\"odinger equation on $u(n)$,  then by the same process, for
$N=4$ (i.e, the fourth-order isospectral flow) with a somewhat long
computation, he has
\begin{eqnarray}
iq_t+q_{xxxx}+4q_{xx}q^*q+2qq^*_{xx}q+4qq^*q_{xx}
+2q_xq^*_xq+6q_xq^*q_x+2qq^*_xq_x+6qq^*qq^*q=0, \label{AKNS01}
\end{eqnarray}
where $$Q=\left(\begin{array}{cc} 0&q\\-q^*&0\end{array}\right)$$ is
an ${\bf m}$-potential of $u(n)$ with $q^*$ being the transposed
conjugate matrix of the complex $k\times(n-k)$-matrix $q$.

\subsection{Modification of generalized bi-Schr\"odinger flows}

We know that the geometric concepts of Schr\"odinger flows and
geometric KdV flows introduced recently in literature (see, for
example, \cite{DW,TeUh,SW,DW1}) are of independent research interest
and enrich the content of research on LIE models. In particular, it
is worthy of note that to make use of the para-K\"ahler structure
provides a new perspective of utilizing the concepts of
Schr\"odinger flows and geometric KdV flows. The results obtained in
\cite{DI,DW1,dinghe} illuminate that the nonlinear Schr\"odinger
equations have the K\"ahler structure rather than the para-K\"ahler
structure, while the nonlinear coupled heat equation and the KdV
equation have the para-K\"ahler structure rather than the K\"ahler
structure. The recognition that the KdV equation possesses the
para-K\"ahler structure seems very appealing. As we shall see below,
the K\"ahler structure on the Grassmannian manifolds $G_{n,k}$ and
$G^k_n$, and the para-K\"ahler structure on the para-Grassmannian
manifold ${\widetilde G}_n^k$ described above will still play
crucial roles in establishing the non-integral third-order model of
the vortex filament on a symmetric Lie algebra.

Let $(M,g)$ be a Riemannian manifold and $(N,h)$ a (pseudo)
Riemannian manifold. Recall that the energy functional and the
bi-energy functional of a map $u: M\to N$ are defined respectively
by
 $$E(u)=\frac{1}{2}\int_M|\nabla u|^2dv_g\quad\quad
\mbox{and}\quad\quad E_2(u)=\frac{1}{2}\int_M\|(d+d^*)^2u\|^2,$$
where $d^*$ denotes the co-differential operator. It is well-known
that harmonic maps (resp. bi-harmonic maps) are defined by critical
points of the energy functional (resp. the bi-energy functional). If
the target manifold $(N,h)$ is a K\"ahler manifold with a compatible
complex structure $J$, then in \cite{DWW1} generalized
bi-Schr\"odinger flow of map $u$ from the Riemannian manifold
$(M,g)$ to the K\"ahler manifold $(N,J,h)$ is defined by the
Hamiltonian gradient flow of the generalized bienergy functional
\begin{eqnarray}
E^*_{\alpha,\beta,\gamma}(u)=\bigg\{\alpha E(u)+\beta E_2(u)+\gamma
\int_M\|\nabla u\|^4\bigg\}.\label{bienergy}
\end{eqnarray} In other words, $u$ satisfies the equation
\begin{eqnarray} u_t=J_u\nabla E^*_{\alpha,\beta,\gamma}(u).
\label{G4}
\end{eqnarray}
Let $(N,K,h)$ be a para-K\"ahler manifold with $K$ as its
para-K\"ahler structure. Similarly, a map $u=u(t,x): [0,T)\times
M\to N$, where $0< T\le\infty$, is called a generalized
bi-Schr\"odinger flow of map from the Riemannian manifold $(M,g)$ to
the para-K\"ahler manifold $(N,K,h)$ if $u$ fulfills the following
para-Hamiltonian gradient flow of the functional
$E^*_{\alpha,\beta,\gamma}(u)$
\begin{eqnarray} u_t=K_u\nabla
E^*_{\alpha,\beta,\gamma}(u). \nn\label{G401}
\end{eqnarray}
We point out that, when $\alpha\not=0$ and $\beta=\gamma=0$,
generalized bi-Schr\"odinger maps reduce to Schr\"odinger flows (see
\cite{DW,TeUh}), or in other words, the equation of Schr\"odinger
flows is just the Hamiltonian gradient flow of the energy functional
$E(u)$.

\medskip

In order for the generalized bi-Schr\"odinger flows to be applicable
in establishing the desired third-order models of the vortex
filament on symmetric Lie algebras, we need to modify the
generalized bi-energy functional (\ref{bienergy}). Our key
observation comes from the following facts. For a map ${u}$ from
${\mathbb R}$ to ${\mathbb S^2}\hookrightarrow {\mathbb R}^{3}$,
${\mathbb H^2}\hookrightarrow {\mathbb R}^{2,1}$ and ${\mathbb
S^{1,1}}\hookrightarrow {\mathbb R}^{2,1}$, we have respectively the
identities:
$$R({u}_x,J_{u}{u}_x)J_{ u}{u}_x=
|{ u}_x|^2{u}_x,\quad\quad R({ u}_x,J_{u}{ u}_x)J_{u}{ u}_x=-|{
u}_x|^2{u}_x$$ and $$R({u}_x,K_{u}{u}_x)K_{u}{u}_x=|{u}_x|^2{ u}_x,
$$ where $R(\cdot,\cdot)$ stands for
the curvature operator on the target manifold $M=\mathbb S^2$,
$\mathbb H^2$ or $\mathbb S^{1,1}$ (refer to \cite{DW1}). Hence, the
third term $\int_M\|\nabla {u}\|^4$ in (\ref{bienergy}) for the
three different cases can be unified by
\begin{eqnarray}
\int_M\langle R({ u}_x,J_{u}{u}_x)J_{u}{u}_x,{
u}_x\rangle,\nn\label{bi*}
\end{eqnarray}
where $\langle\cdot,\cdot\rangle$ and $J_u$ denote respectively the
inner product induced from the standard (pseudo) Riemannian metric
and the complex (or para-complex) structure) on $M={\mathbb S^2}$,
${\mathbb H^2}$ or ${\mathbb S^{1,1}}$.

Based on this observation, we see that the generalized bi-energy
functional (\ref{bienergy}) of smooth maps $u$ from a Riemannian
manifold $(M,g)$ to a K\"ahler manifold
$(N,J,h)$ might be modified as follows. Let $\{e_1,e_2,\cdots,e_n\}$
be a local frame of $(M,g)$, and the metric $g$ in the frame is expressed
by $g=(g_{ij})$ and $(g^{ij})$ its inverse. Then we define
\begin{eqnarray}E_{\alpha,\beta,\gamma}(u)=\bigg\{\alpha E(u)+\beta
E_2(u)+\gamma \int_M\sum_{i,j,k,l=1}^ng^{ij}g^{kl}\langle
R({\nabla_{e_i} u}, J_{ u}{\nabla_{e_j} u})J_{u}{\nabla_{e_k}
u},{\nabla_{e_l} u}\rangle\bigg\},\label{energynew}\end{eqnarray}
where $R(\cdot,\cdot)$ is the curvature operator on the target
manifold $N$ and $\nabla_{e_k}$ the covariant derivative on the
pull-back bundle $u^{-1}TN$ induced from the Levi-Civita connection
on $N$. We should point out that the density function
 $$\sum\limits_{i,j,k,l=1}^ng^{ij}g^{kl}\langle
 R({\nabla_{e_i} u}, J_{ u}{\nabla_{e_j} u})J_{u}{\nabla_{e_k}u},{\nabla_{e_l} u}\rangle$$
of the $\gamma$-term in (\ref{energynew}) is independent of the
choice of a local frame and hence the definition of the generalized
bi-energy function (\ref{energynew}) is of significance.

If the target manifold is a para-K\"ahler manifold denoted by $(N,K,h)$,
similarly we should define the corresponding functional as follows
\begin{eqnarray*}E_{\alpha,\beta,\gamma}(u)=\bigg\{\alpha E(u)+\beta
E_2(u)+\gamma \int_M\sum_{i,j,k,l=1}^ng^{ij}g^{kl}\langle
R({\nabla_{e_i} u}, K_{ u}{\nabla_{e_j} u})K_{u}{\nabla_{e_k}
u},{\nabla_{e_l} u}\rangle\bigg\}.\label{energynew*}\end{eqnarray*}

\begin{Definition} Let $(N,J,h)$ (resp. $(N,K,h)$)  be a
K\"ahler (resp. para-K\"ahler) manifold with $J$ (resp. $K$) as its
K\"ahler (resp. para-K\"ahler) structure. A map $u=u(t,x):
[0,T)\times M\to N$, where $0< T\le\infty$, is called a {\bf
generalized bi-Schr\"odinger flow} from the Riemannian manifold
$(M,g)$ to $(N,J,h)$ (resp. $(N,K,h)$) if $u$ satisfies the following
Hamiltonian gradient flow
\begin{eqnarray} u_t=J_u\nabla
E_{\alpha,\beta,\gamma}(u)\quad\quad \left(\hbox{resp.}~
u_t=K_u\nabla E_{\alpha,\beta,\gamma}(u)\right), \label{G403}
\end{eqnarray}where $E_{\alpha,\beta,\gamma}(u)$ is the
generalized bi-energy functional defined by (\ref{energynew}).

A map $u=u(t,x): [0,T)\times M\to N$ is called a {\bf
bi-Schr\"odinger flow} from the Riemannian manifold $(M,g)$ to
$(N,J,h)$ if $u$ satisfies the following Hamiltonian gradient flow
\begin{eqnarray} u_t=J_u\nabla
E_{0,1,\gamma}(u)\quad\quad(\hbox{resp.}~ u_t=K_u\nabla
E_{0,1,\gamma}(u)). \label{Gwang}
\end{eqnarray}
\end{Definition}

We remark that without the contribution of the $\gamma$-term
appearing in the generalized bi-energy functional
$E_{\alpha,\beta,\gamma}(u)$, one cannot obtain completely the
desired model on a symmetric Lie algebra, as we shall see in \S3.
This gives the reason why the $\gamma$-term which relies seriously
on the intrinsic geometry of the target manifold in the functional
$E_{\alpha,\beta,\gamma}(u)$ is indispensable. However, one also
notes that the concept of generalized bi-Schr\"odinger maps is
mainly based on the bi-energy and hence bi-Schr\"odinger maps.
Therefore, the bi-energy plays an essential role in this respect.
Furthermore, comparing $E_{\alpha,\beta,\gamma}(u)$
(\ref{energynew}) with $E^*_{\alpha,\beta,\gamma}(u)$
(\ref{bienergy}) in the case that the start manifold is $\mathbb R$
and the target manifold is to be $\mathbb H^2$, there is an
additional sign in the third term of the generalized bi-energy
functional $E_{\alpha,\beta,\gamma}(u)$. The so-called mysterious
question  why the parameters $\beta$ and $\gamma$ have to satisfy
the different conditions arises from the geometric explanation of
the second to fourth order integrable equations in the AKNS
hierarchy by generalized bi-Schr\"odinger flows in \cite{DLW}
becomes clear now. By using the present concept of generalized
bi-Schr\"odinger flows, the parameters $\beta$ and $\gamma$ do
indeed satisfy the same condition: $\gamma=-\frac{\beta}{8}$ (one
may refer to \S3 of \cite{DLW} for details). Hence, the modification
looks very natural and effective. In the following context,
generalized bi-Schr\"odinger flows are always referred to in the
present sense, unless otherwise specified.

\medskip

{\bf Examples}:

a) The equation of generalized bi-Schr\"odinger flows from $\mathbb
R$ to ${\mathbb S^2}\hookrightarrow {\mathbb R}^{3}$ is just the one
displayed in \cite{DWW1,DLW}, i.e. for ${\bf s}=(s_1,s_2,s_3)\in
\mathbb R^3$, \begin{eqnarray}{\bf s}_t={\bf
s}\times\left(-\alpha{\bf s}_{xx}+\beta {\bf
s}_{xxxx}-(4\gamma-2\beta)(|{\bf s}_x|^2 {\bf s}_x)_x\right),\quad
|{\bf s}|^2:=s_1^2+s_2^2+s_3^2=1.\label{S2}\end{eqnarray}

b) The equation  of generalized bi-Schr\"odinger flows from $\mathbb
R$ to ${\mathbb H^2}\hookrightarrow {\mathbb R}^{2,1}$ is the one
displayed in \cite{DWW1,DLW} with $\gamma$ changed by $-\gamma$,
i.e. for ${\bf s}=(s_1,s_2,s_3)\in \mathbb R^{2,1}$,
\begin{eqnarray}
\textbf{s}_t=\textbf{s}\dot{\times}\left(-\alpha\textbf{s}_{xx}+\beta\textbf{s}_{xxxx}
+(4\gamma-2\beta)(|\textbf{s}_x|^2\textbf{s}_x)_x\right),\ \ \
|\textbf{s}|^2:=s_1^2+s_2^2-s_3^2=-1.\label{H2}
\end{eqnarray}

c) The equation  of generalized bi-Schr\"odinger flows from $\mathbb
R$ to ${\mathbb S^{1,1}}\hookrightarrow {\mathbb R}^{2,1}$ is the
same one as displayed in \cite{DLW}, i.e. \begin{eqnarray}{\bf
s}_t={\bf s}\dot{\times}\left(-\alpha{\bf s}_{xx}+\beta{\bf
s}_{xxxx}-(4\gamma-2\beta)(| {\bf s}_x|^2 {\bf s}_x)_x\right),\quad
|{\bf s}|^2:=s_1^2+s_2^2-s_3^2=1.\label{S11}\end{eqnarray}

It is proved in \cite{DWW1,DLW} that the equations
(\ref{S2}), (\ref{H2}) and (\ref{S11}) of the generalized
bi-Schr\"odinger flows are respectively equivalent to the
non-integrable Fukumoto-Moffatt's model (\ref{I0}) in $\mathbb R^{3}$, the timelike
third-order model (\ref{I1}) and the spacelike third-order model
(\ref{I2}) of the vortex filament in the Minkowski space $\mathbb R^{2,1}$.

\section{Generalized bi-Schr\"odinger flows into symmetric spaces}

In this section, we shall use our new generalized bi-Schr\"odinger
flows introduced above to deduce the equations of the generalized
bi-Schr\"odinger flows to the compact K\"ahler Grassmannian
manifolds $G_{n,k}$, the noncompact K\"ahler Grassmannian manifolds
$G_{n}^k$ and the para-K\"ahler Grassmannian manifolds ${\widetilde
G}_{n,k}$ respectively. These models return respectively to the
non-integrable Fukumoto-Moffatt's model (\ref{I0}) in $\mathbb R^3$,
the timelike third-order model (\ref{I1}) and the spacelike
third-order model (\ref{I2}) of the vortex filament in $\mathbb
R^{2,1}$ when the target manifold $M$ goes back to ${\mathbb S^2}$,
${\mathbb H^2}$ and ${\mathbb S^{1,1}}$ correspondingly. Therefore,
they are naturally the Fukumoto-Moffatt's counterparts of the vortex
filament on symmetric Lie algebras.

First of all, we note that $\langle A,B\rangle=-\hbox{tr}(AB)$ is
the bi-invariant inner product on $u(n)$ (refer to, for example,
\cite{TeUh}), $\langle A,B\rangle=\hbox{tr}(AB)$ is the bi-invariant
inner product on $u(k,n-k)$ (refer to, for example, \cite{He}) and
$\langle A,B\rangle=\hbox{tr}(AB)$ is the bi-invariant pseudo inner
product on $gl(n,\mathbb R)$. Now for a map $u$ from $\mathbb R$ to
the K\"ahler Grassmannian manifolds $G_{n,k}=\left\{E^{-1}\sigma_3
E~|~ \forall ~E\in U(n)\right\}$ of compact type, or to the K\"ahler
Grassmannian manifolds $G_{n}^k=\left\{E^{-1}\sigma_3 E~|~ \forall
~E\in U(k,n-k)\right\}$ of noncompact type and to the para-K\"ahler
Grassmannian manifolds ${\widetilde G}_{n,k}=\left\{E^{-1}\sigma_3
E~|~ \forall ~E \in\right.$ $\left.GL(n,\mathbb R)\right\}$, we see
that the $\gamma$-term of the generalized bi-energy functional
$E_{\alpha,\beta,\gamma}(u)$ given by (\ref{energynew}) becomes
\begin{eqnarray}
\int_{\mathbb R} \langle R({ u}_x,J_{ u}{ u}_x)J_{ u}{ u}_x,{
u}_x\rangle ~~~~\left(\hbox{resp}~~\int_{\mathbb R}\langle R({
u}_x,K_{ u}{ u}_x)K_{ u}{ u}_x,{ u}_x\rangle\right),\nn
\end{eqnarray}
where $x$ stands for the standard coordinate on $\mathbb R$ and
$u_x=\nabla_{\frac{\partial}{\partial x}}u$. Now we have

\begin{Theorem} \label{Thm1}
The equation of the generalized bi-Schr\"{o}dinger flows (\ref{G403})
from $\mathbb R$ into the symmetric space $M$ reads
\begin{eqnarray}
\varphi_t =\bigg[\varphi,\,\, -\alpha\varphi_{xx}
+\beta\varphi_{xxxx}+(4\gamma-2\beta)\left(\varphi_x\varphi^{-1}\varphi_x
\varphi^{-1}\varphi_x\right)_x\bigg],\label{bi0}
\end{eqnarray}
where $M$ is one of $G_{n,k}$, $G_{n}^k$ or ${\widetilde G}_{n,k}$,
this corresponds respectively to the Hermitian symmetric Lie algebra
$u(n)$ with index $k$ of compact type, the Hermitian symmetric Lie
algebra $u(k,n-k)$ with index $k$ of noncompact type and the
para-Hermitian symmetric Lie algebra $gl(n,\mathbb R)$ with index
$k$.
\end{Theorem}
{\bf Proof}: We will reach the conclusion case by case as follows.

1) Let $M=G_{n,k}=\left\{E^{-1}\sigma_3 E~|~ \forall ~E\in
U(n)\right\}$, the Grassmannian manifolds of compact type, and
$\varphi=E^{-1}(t,x)\sigma_3E(t,x)$ be a map from the line $\mathbb
R$ to $M=\left\{E^{-1}\sigma_3 E~|~ \forall ~E\in U(n)\right\}$,
where $\sigma_3$ is given by (\ref{1}). Without loss of
generality, we may assume that $E$ satisfies: $E_x=PE$ for some
$P\in {\bf m}$, where ${\bf m}$ fits the symmetric condition:
$u(n)={\bf k}\oplus {\bf m}$. In fact, if $E$ does not meet the
requirement, we may make a transform: $$E\to
\left(\begin{array}{cc}A&0\\0&B\end{array}\right)E$$ (this is because
the form $E^{-1}\sigma_3E$ is invariant up to the transform) such
that by suitably choosing $A$ and $B$ (through solving a linear
differential system of $A$ and $B$)  $P$ can be modified so that
the new $P$ satisfies $P\in {\bf m}$. Based on this fact, we have
\begin{eqnarray}
\varphi^{-1}&=&4\varphi^*=-4\varphi,~~~\Longrightarrow
~(\varphi^{-1})^2=-4I,\label{II16}\\
\varphi_x&=&E^{-1}[\sigma_3,P]E, ~~\Longrightarrow
~[\varphi,\varphi_x]=-E^{-1}PE,\label{II17}\\
\varphi_x\varphi^{-1}&=&-2E^{-1}PE,\quad\quad\varphi^{-1}\varphi_x=2E^{-1}PE.\label{II18}
\end{eqnarray}
Combining (\ref{II17}) with (\ref{II18}), we have
\begin{eqnarray}
[\varphi,\,\varphi_x]=J_{\varphi}\varphi_x=\frac{1}{2}
\varphi_x\varphi^{-1}\label{II19}
\end{eqnarray}
and
\begin{eqnarray}
\varphi^2_x=E^{-1}[\sigma_3,P]^2E=E^{-1}P^2E=\frac{1}{4}
\varphi_x\varphi^{-1}\varphi_x\varphi^{-1}.\label{II20}
\end{eqnarray}

Now we rewrite the generalized bi-energy $E_{\alpha,\beta,\gamma}(\varphi)$ as
\begin{eqnarray}
E_{\alpha,\beta,\gamma}(u)&=&\bigg\{\alpha E(u)+\beta E_2(u)+\gamma
\int_{\mathbb R} \langle R({ u}_x,J_{ u}{ u}_x)J_{ u}{ u}_x,{
u}_x\rangle\bigg\}\nn\\
&:=&\alpha E(\varphi)+\beta E_2(\varphi)+\gamma {\widetilde
E}(\varphi).\nn
\end{eqnarray}
Hence the calculation of the gradient of the generalized bi-energy
$E_{\alpha,\beta,\gamma}(\varphi)$  reduces to calculating the
gradient of $E(\varphi)$, $E_2(\varphi)$ and ${\widetilde
E}(\varphi)$ respectively. Firstly, we come to calculate explicitly
the gradient of the energy functional $E(\varphi)$. It is well known
(see, for example, \cite{TeUh}) that
\begin{eqnarray}
\nabla E=\tau(\varphi)=-\nabla_{\varphi_x}\varphi_x, \label{II12}
\end{eqnarray} where $\tau(\varphi)$ stands for the tension field of
$\varphi$. One may verify directly by using (\ref{II17}) and
(\ref{II18}) that
\begin{eqnarray}\nabla_{\varphi_x}\varphi_x&=&\varphi_{xx}+E^{-1}[P,\,[\sigma_3,\,P]]E
~~~\hbox{(taking the tangent
part)}\nn\\
&=&\varphi_{xx}-\varphi_x\varphi^{-1}\varphi_x.\label{II121}\end{eqnarray}

Next, by using the identities (\ref{II18},\ref{II16}) and the fact
(\ref{II121}) (with (\ref{II12})), we see that the bi-energy
functional $E_2(\varphi)$ can be re-written as
\begin{eqnarray}
E_2(\varphi) &=&\frac{1}{2}\int_{\mathbb
R}\|(d+d^*)^2\varphi\|^2=\frac{1}{2}\int_{\mathbb
R}\langle\tau(\varphi),\,\tau(\varphi)\rangle\nn\\
&=&\frac{1}{2}\int_{\mathbb R}\bigg(
\langle\varphi_{xx},\varphi_{xx}\rangle-2\langle\varphi_{xx},\varphi_x\varphi^{-1}\varphi_x\rangle+
\langle\varphi_x\varphi^{-1}\varphi_x,\varphi_x\varphi^{-1}\varphi_x\rangle\bigg)\nonumber\\
&=&\frac{1}{2}\int_{\mathbb
R}\langle\varphi_{xx},\varphi_{xx}\rangle-\int_{\mathbb
R}\langle\varphi_{xx},\varphi_x\varphi^{-1}\varphi_x\rangle+\frac{1}{2}\int_{\mathbb
R}\langle\varphi_x\varphi^{-1}\varphi_x,\varphi_x\varphi^{-1}\varphi_x\rangle\nonumber\\
&:=&E_{21}(\varphi)-E_{22}(\varphi)+E_{23}(\varphi),\label{II13}
\end{eqnarray}
where \begin{eqnarray} E_{21}(\varphi)&=&\frac{1}{2}\int_{\mathbb
R}\langle\varphi_{xx},\,\varphi_{xx}\rangle,\nn\\E_{22}(\varphi)&=&\int_{\mathbb
R}\langle\varphi_{xx},\,\varphi_x\varphi^{-1}\varphi_x\rangle,\nn\\E_{23}(\varphi)&=&\frac{1}{2}\int_{\mathbb
R}
\langle\varphi_x\varphi^{-1}\varphi_x,\,\varphi_x\varphi^{-1}\varphi_x\rangle.\nn\end{eqnarray}
The calculation of the gradient of the bi-energy functional
$E_2(\varphi)$ comes down to calculating the gradients of
$E_{21}(\varphi)$, $E_{22}(\varphi)$ and $E_{23}(\varphi)$
respectively.

It is easy to see that the variation of the functional
$E_{21}(\varphi)$ is
\begin{eqnarray}
\lim_{\varepsilon\to 0}\frac{E_{21}(\varphi+\varepsilon
h)-E_{21}(\varphi)}{\varepsilon} =~\int_{\mathbb R}\langle\varphi_{xx},
h_{xx}\rangle=\int_{\mathbb R}\langle\varphi_{xxxx},h\rangle.\nn
\end{eqnarray}
This indicates that
\begin{eqnarray}
\nabla E_{21}&=& \varphi_{xxxx}+\hat {\bf k}_1\label{II14}
\end{eqnarray}
for some ${\bf k}$-valued smooth function $\hat {\bf k}_1$. To
calculate the gradient of $E_{22}(\varphi)$, we reexpress it as
follows,
\begin{eqnarray}
E_{22}(\varphi)&=&\int_{\mathbb
R}\langle\varphi_{xx},\varphi_x\varphi^{-1}\varphi_x\rangle\nn\\
&=&-\int_{\mathbb
R}\langle\varphi_{xx},\varphi^{-1}\varphi_x\varphi_x\rangle\nn\\
&=&\frac{1}{4}\int_{\mathbb
R}\hbox{tr}\left(\varphi_{xx}\varphi^{-1}\varphi_x\varphi^{-1}\varphi_x\varphi^{-1}\right).\label{II122}
\end{eqnarray}
Here we have used the formulae (\ref{II18}) and (\ref{II20}) as well as the
invariant inner product on $u(n)$ mentioned above. Therefore
\begin{eqnarray}
&&\lim_{\varepsilon\to 0}\frac{E_{22}(\varphi+\varepsilon
h)-E_{22}(\varphi)}{\varepsilon}\nonumber\\
&=&\frac{1}{4}\int_{\mathbb
R}\hbox{tr}\bigg(h_{xx}\varphi^{-1}\varphi_x\varphi^{-1}\varphi_x\varphi^{-1}+
\varphi_{xx}\varphi^{-1}\varphi_x\varphi^{-1}h_x\varphi^{-1}+
\varphi_{xx}\varphi^{-1}h_x\varphi^{-1}\varphi_x\varphi^{-1}
\nn\\
&&-\left(\varphi^{-1}\varphi_{x}\varphi^{-1}\varphi_x\varphi^{-1}\varphi_{xx}\varphi^{-1}
+\varphi^{-1}\varphi_{xx}\varphi^{-1}\varphi_x\varphi^{-1}\varphi_x\varphi^{-1}
+\varphi^{-1}\varphi_{x}\varphi^{-1}\varphi_{xx}\varphi^{-1}\varphi_x\varphi^{-1}
\right)h\bigg)\nn\\
&=&\frac{1}{4}\int_{\mathbb
R}\hbox{tr}\bigg(\bigg[\left(\varphi^{-1}\varphi_x\varphi^{-1}\varphi_x\varphi^{-1}\right)_{xx}-
\left(\varphi^{-1}\varphi_{xx}\varphi^{-1}\varphi_x\varphi^{-1}\right)_x-
\left(\varphi^{-1}\varphi_x\varphi^{-1}\varphi_{xx}\varphi^{-1}\right)_x\nn\\
&&~~~ -
\left(\varphi^{-1}\varphi_x\varphi^{-1}\varphi_x\varphi^{-1}\varphi_x\varphi^{-1}\right)_x\bigg]h\bigg).
\nn
\end{eqnarray}
Hence we obtain
\begin{eqnarray}
\nabla
E_{22}&=&-\frac{1}{4}\bigg(\left(\varphi^{-1}\varphi_x\varphi^{-1}\varphi_x\varphi^{-1}\right)_{xx}-
\left(\varphi^{-1}\varphi_{xx}\varphi^{-1}\varphi_x\varphi^{-1}\right)_x-\nn\\
&&~~~\left(\varphi^{-1}\varphi_x\varphi^{-1}\varphi_{xx}\varphi^{-1}\right)_x
-
\left(\varphi^{-1}\varphi_x\varphi^{-1}\varphi_x\varphi^{-1}\varphi_x\varphi^{-1}\right)_x\bigg)+\hat{\bf
k}_2\nn\\
&=&
\left(\varphi^{-1}\varphi_x\varphi^{-1}\varphi_x\varphi^{-1}\varphi_x\varphi^{-1}\right)_x+\hat{\bf
k}_2\label{II15}
\end{eqnarray}
for some  ${\bf k}$-valued smooth function $\hat {\bf k}_2$.

Let's now treat $E_{23}(\varphi)$. One sees that
\begin{eqnarray}
E_{23}(\varphi)&=&\frac{1}{2}\int_{\mathbb R}
\langle\varphi_x\varphi^{-1}\varphi_x,\varphi_x\varphi^{-1}\varphi_x\rangle\nn\\
&=&-\frac{1}{2}\int_{\mathbb R}\hbox{tr}\bigg(
\varphi_x\varphi^{-1}\varphi_x\varphi_x\varphi^{-1}\varphi_x\bigg)\nn\\
&=&\frac{1}{8}\int_{\mathbb
R}\hbox{tr}\bigg(\varphi_x\varphi^{-1}\varphi_x\varphi^{-1}\varphi_x\varphi^{-1}\varphi_x\varphi^{-1}\bigg).
\label{E23}
\end{eqnarray}
Here we have used the formula (\ref{II20}). Before we calculate the
gradient of $E_{23}(\varphi)$, let's treat the $\gamma$-term
$$
{\widetilde E}(\varphi)=\int_{\mathbb
R}\langle R(\varphi_x,J_{\varphi}\varphi_x)J_{\varphi}\varphi_x,\varphi_x\rangle$$
 of the generalized bi-energy functional
$E_{\alpha,\beta,\gamma}(\varphi)$. We know that (refer to
\cite{ppetersen})
$$R(\varphi_x,J_{\varphi}\varphi_x)J_{\varphi}\varphi_x=
[J_{\varphi}\varphi_x,[\varphi_x,J_{\varphi}\varphi_x]].$$ Hence, by
using (\ref{II16},\,\ref{II19}) and (\ref{II20}), we obtain
\begin{eqnarray}
R(\varphi_x,J_{\varphi}\varphi_x)J_{\varphi}\varphi_x&=&
\frac{1}{4}[\varphi_x\varphi^{-1},[\varphi_x,\varphi_x\varphi^{-1}]]=
\frac{1}{4}[\varphi_x\varphi^{-1},\varphi_x^2\varphi^{-1}-\varphi_x\varphi^{-1}\varphi_x]
\nn\\
&=&\frac{1}{4}\bigg(2\varphi_x\varphi^{-1}\varphi_x^2\varphi^{-1}-\varphi_x\varphi^{-1}\varphi_x\varphi^{-1}
\varphi_x-\varphi^2_x\varphi^{-1}\varphi_x\varphi^{-1}\bigg)\nn\\
&=&\frac{1}{4}\bigg(2\varphi_x\varphi^{-1}\left(\frac{1}{4}
\varphi_x\varphi^{-1}\varphi_x\varphi^{-1}\right)\varphi^{-1}-\varphi_x\varphi^{-1}\varphi_x\varphi^{-1}
\varphi_x\nn\\
&&~~~~~~-\left(\frac{1}{4}
\varphi_x\varphi^{-1}\varphi_x\varphi^{-1}\right)\varphi^{-1}\varphi_x\varphi^{-1}\bigg)\nn\\
&=&-\varphi_x\varphi^{-1}\varphi_x\varphi^{-1} \varphi_x\nn
\end{eqnarray}
and
\begin{eqnarray}
{\widetilde E}(\varphi)&=&\int_{\mathbb
R}\langle R(\varphi_x,J_{\varphi}\varphi_x)J_{\varphi}\varphi_x,\varphi_x\rangle\nn\\
&=&-\int_{\mathbb
R}\hbox{tr}\bigg(\left(-\varphi_x\varphi^{-1}\varphi_x\varphi^{-1}
\varphi_x\right)\cdot \varphi_x\bigg)\nn\\
&=&\int_{\mathbb
R}\hbox{tr}\bigg(\left(-\varphi_x\varphi^{-1}\varphi_x\varphi^{-1}
\varphi_x\right)\cdot \left(\frac{1}{4}\varphi^{-1}\right)_x\bigg)\nn\\
&=&\frac{1}{4}\int_{\mathbb
R}\hbox{tr}\bigg(\varphi_x\varphi^{-1}\varphi_x\varphi^{-1}\varphi_x\varphi^{-1}\varphi_x\varphi^{-1}\bigg)
\nn\\
&=&{2}E_{23}(\varphi).\label{II21}
\end{eqnarray}
It is direct to see that the variation of ${\widetilde E}(\varphi)$
is
\begin{eqnarray}
&&\lim_{\varepsilon\to 0}\frac{{\widetilde E}(\varphi+\varepsilon
h)-{\widetilde E}(\varphi)}{\varepsilon}\nn\\
&=&\frac{1}{4}\lim_{\varepsilon\to 0}\int_{\mathbb
R}\frac{\hbox{tr}\bigg(\left((\varphi+\varepsilon
h)_x(\varphi+\varepsilon h)^{-1}\right)^4-
\left(\varphi_x\varphi^{-1}\right)^4\bigg)}{\varepsilon}\nn\\
&=&\frac{1}{4}\int_{\mathbb R}\hbox{tr}\bigg(
\{\varphi_x\varphi^{-1}\varphi_x\varphi^{-1},
\{\varphi_x\varphi^{-1},h_x\varphi^{-1}-\varphi_x\varphi^{-1}h\varphi^{-1}\}\}\bigg)\nn\\
&=&\frac{1}{4} \int_{\mathbb
R}\hbox{tr}\bigg(-4(\varphi_x\varphi^{-1}\varphi_x\varphi^{-1}\varphi_x\varphi^{-1})_xh\varphi^{-1}-
4(\varphi_x\varphi^{-1})^4h\varphi^{-1}\bigg)\nn\\
&=&\int_{\mathbb
R}\langle(\varphi_x\varphi^{-1}\varphi_x\varphi^{-1}\varphi_x\varphi^{-1})_x+
(\varphi_x\varphi^{-1})^4,\, h\varphi^{-1}\rangle.\label{II22}
\end{eqnarray}
Hence
\begin{eqnarray}
\nabla {\widetilde
E}=\varphi^{-1}\bigg((\varphi_x\varphi^{-1}\varphi_x\varphi^{-1}\varphi_x\varphi^{-1})_x+
(\varphi_x\varphi^{-1})^4+\hat {\bf k}_3\bigg)\label{II23}
\end{eqnarray}
and
\begin{eqnarray}
\nabla
{E}_{23}=\frac{1}{2}\bigg((\varphi^{-1}\varphi_x\varphi^{-1}\varphi_x\varphi^{-1}\varphi_x\varphi^{-1})_x+
\hat{\bf k}_3\bigg)\label{II230}
\end{eqnarray}
for some  ${\bf k}$-valued smooth function $\hat {\bf k}_3$.

Combining (\ref{II14},\ref{II15}) with (\ref{II23}) and
(\ref{II230}), we see that the equation of the generalized
bi-Schr\"odinger flow of maps from ${\mathbb R}$ to the K\"ahler
Grassmannian manifolds $G_{n,k}$ of compact type is
\begin{eqnarray}
\varphi_t&=&\bigg[\varphi,\, \alpha\nabla E+\beta\nabla E_2+\gamma\nabla {\widetilde E}\bigg]\nn\\
&=&\bigg[\varphi,\, \alpha\nabla E+\beta\nabla
\left(E_{21}-E_{22}\right)+
\left(\beta+2\gamma\right) \nabla E_{23}\bigg]\nn\\
&=&\bigg[\varphi,\, -{\alpha}\left(\varphi_{xx}-\varphi_x\varphi^{-1}\varphi_x\right)
+\beta\left(\varphi_{xxxx}+\hat{\bf
k}_1-(\varphi^{-1}\varphi_x\varphi^{-1}\varphi_x\varphi^{-1}\varphi_x\varphi^{-1})_x
+\hat{\bf k}_2\right)\nn\\
&&~~~
+\frac{1}{2}(\beta+2\gamma)\left((\varphi^{-1}\varphi_x\varphi^{-1}\varphi_x\varphi^{-1}\varphi_x\varphi^{-1})_x
+\hat{\bf
k}_3\right)\bigg]\nn\\
&=&\bigg[\varphi,\,-\alpha\varphi_{xx}
+\beta\varphi_{xxxx}+\left(\gamma-\frac{\beta}{2}\right)
\varphi^{-1}\left(\varphi_x\varphi^{-1}\varphi_x\varphi^{-1}\varphi_x\right)_x\varphi^{-1}\bigg]
\nn\\
&=&\bigg[\varphi,\, -\alpha\varphi_{xx}
+\beta\varphi_{xxxx}+(4\gamma-2\beta)\left(\varphi_x\varphi^{-1}\varphi_x\varphi^{-1}\varphi_x\right)_x\bigg],\nn
\end{eqnarray}
which is just Eq.(\ref{bi0}). Here we have used the identity
(\ref{II16}) and the symmetric relation: $[{\bf k},{\bf k}]=0$. This
completes the proof of the theorem in the case $M=G_{n,k}$.

\medskip

2) For $M=G_{n}^k=\left\{E^{-1}\sigma_3 E~|~ \forall ~E\in
U(k,n-k)\right\}$, i.e. the K\"ahler Grassmannian manifolds of
noncompact type, one may verify that the formulas (\ref{II16}~$\sim$
\ref{II20}) are still valid in the present case. One notes that the
metric used for the K\"ahler Grassmannian manifolds $G_{n}^k$ of
noncompact type, i.e. in the case of the symmetric Lie algebra
$u(k,n-k)$, is given by $\langle A,B\rangle = \hbox{tr}(AB)$,
$\forall A,B\in u(k,n-k)$. Following the line of proof of Theorem
\ref{Thm1}, we may similarly have
\begin{eqnarray}
\nabla E=\tau(\varphi)=-\nabla_{\varphi_x}\varphi_x=-\varphi_{xx}+
\varphi_x\varphi^{-1}\varphi_x,\nn\label{II121non}\end{eqnarray}
\begin{eqnarray}
E_{22}(\varphi)=-\frac{1}{4}\int_{\mathbb
R}\hbox{tr}\left(\varphi_{xx}\varphi^{-1}\varphi_x\varphi^{-1}\varphi_x\varphi^{-1}\right)\nn\label{II122non}
\end{eqnarray}
and \begin{eqnarray} {\widetilde E}(\varphi)&=&\int_{\mathbb
R}\langle R(\varphi_x,J_{\varphi}\varphi_x)J_{\varphi}\varphi_x,\, \varphi_x\rangle\nn\\
&=&-\frac{1}{4}\int_{\mathbb
R}\hbox{tr}\bigg(\varphi_x\varphi^{-1}\varphi_x\varphi^{-1}\varphi_x\varphi^{-1}\varphi_x\varphi^{-1}\bigg)
\nn\\
&=&{2}E_{23}(\varphi).\nn\label{II21no}
\end{eqnarray}
Here we have used the invariant (pseudo) inner product on
$u(k,n-k)$. Therefore,
\begin{eqnarray}
\nabla E_{21}&=& \varphi_{xxxx}+\hat {\bf k}_1\label{II14no}
\end{eqnarray}
for some ${\bf k}$-valued smooth function $\hat {\bf k}_1$,
\begin{eqnarray}
\nabla E_{22}=
\left(\varphi^{-1}\varphi_x\varphi^{-1}\varphi_x\varphi^{-1}\varphi_x\varphi^{-1}\right)_x+\hat{\bf
k}_2\label{II15no}
\end{eqnarray}
for some  ${\bf k}$-valued smooth function $\hat {\bf k}_2$ and
\begin{eqnarray}
\nabla
{E}_{23}=\frac{1}{2}\bigg((\varphi^{-1}\varphi_x\varphi^{-1}\varphi_x\varphi^{-1}\varphi_x\varphi^{-1})_x+
\hat{\bf k}_3\bigg)\label{II230no}
\end{eqnarray}
for some  ${\bf k}$-valued smooth function $\hat {\bf k}_3$. By
using the formulae (\ref{II14no},\ref{II15no}) and (\ref{II230no}),
we see that the equation of the generalized bi-Schr\"odinger flow of
maps from ${\mathbb R}$ to the K\"ahler Grassmannian manifolds
$G^k_{n}$ of noncompact type is
\begin{eqnarray}
\varphi_t&=&\bigg[\varphi,\alpha\nabla E+\beta\nabla E_2+\gamma\nabla {\widetilde E}\bigg]\nn\\
&=&\bigg[\varphi,\, \alpha\nabla E+\beta\nabla
\left(E_{21}-E_{22}\right)+
\left(\beta+2\gamma\right) \nabla E_{23}\bigg]\nn\\
&=&\bigg[\varphi,-{\alpha}\left(\varphi_{xx}-\varphi_x\varphi^{-1}\varphi_x\right)
+\beta\left(\varphi_{xxxx}+\hat{\bf
k}_1-(\varphi^{-1}\varphi_x\varphi^{-1}\varphi_x\varphi^{-1}\varphi_x\varphi^{-1})_x
+\hat{\bf k}_2\right)\nn\\
&&~~~
+\frac{1}{2}(\beta+2\gamma)\left((\varphi^{-1}\varphi_x\varphi^{-1}\varphi_x\varphi^{-1}\varphi_x\varphi^{-1})_x
+\hat{\bf
k}_3\right)\bigg]\nn\\
&=&\bigg[\varphi,-\alpha\varphi_{xx}
+\beta\varphi_{xxxx}+\left(\gamma-\frac{\beta}{2}\right)
\varphi^{-1}\left(\varphi_x\varphi^{-1}\varphi_x\varphi^{-1}\varphi_x\right)_x\varphi^{-1}\bigg]
\nn\\
&=&\bigg[\varphi,-\alpha\varphi_{xx}
+\beta\varphi_{xxxx}+(4\gamma-2\beta)\left(\varphi_x\varphi^{-1}\varphi_x\varphi^{-1}\varphi_x\right)_x\bigg],\nn
\end{eqnarray}
which also has the same form as Eq.(\ref{bi0}).

3) For the case $M={\widetilde G}_{n,k}=\left\{E^{-1}\sigma_3 E~|~
\forall ~E\in GL(n,\mathbb R)\right\}$, i.e. the para-K\"ahler
Grassmannian manifolds, using the invariant metric on $gl(n)$ or
$sl(2n)$, one may verify that the formulas (\ref{II16}-\ref{II20})
become
\begin{eqnarray}
\varphi^{-1}&=&4\varphi,~~~\Longrightarrow
~(\varphi^{-1})^2=4I,\label{II16p}\\
\varphi_x&=&E^{-1}[\sigma_3,P]E, ~~\Longrightarrow
~[\varphi,\varphi_x]=E^{-1}PE,\label{II17p}\\
\varphi_x\varphi^{-1}&=&-2E^{-1}PE,~\varphi^{-1}\varphi_x=2E^{-1}PE.\label{II18p}
\end{eqnarray}
Combining (\ref{II17p}) with (\ref{II18p}), we have
\begin{eqnarray}
[\varphi,\varphi_x]=J_{\varphi}\varphi_x=-\frac{1}{2}
\varphi_x\varphi^{-1}\label{II19p}
\end{eqnarray}
and
\begin{eqnarray}
\varphi^2_x=E^{-1}[\sigma_3,P]^2E=-E^{-1}P^2E=-\frac{1}{4}
\varphi_x\varphi^{-1}\varphi_x\varphi^{-1}.\label{II20p}
\end{eqnarray}
As in the case of $u(n)$, we similarly have that
\begin{eqnarray}
\nabla E&=&\tau(\varphi)=-\nabla_{\varphi_x}\varphi_x=-\varphi_{xx}+
\varphi_x\varphi^{-1}\varphi_x,\nn\label{II12para}\end{eqnarray}
\begin{eqnarray}
E_{22}(\varphi)=-\frac{1}{4}\int_{\mathbb
R}\hbox{tr}\left(\varphi_{xx}\varphi^{-1}\varphi_x\varphi^{-1}\varphi_x\varphi^{-1}\right),
\nn\label{II122para}
\end{eqnarray}
and \begin{eqnarray} {\widetilde E}(\varphi)&=&\int_{\mathbb
R}\langle R(\varphi_x,J_{\varphi}\varphi_x)J_{\varphi}\varphi_x,\,\varphi_x\rangle\nn\\
&=&\frac{1}{4}\int_{\mathbb
R}\hbox{tr}\bigg(\varphi_x\varphi^{-1}\varphi_x\varphi^{-1}\varphi_x\varphi^{-1}\varphi_x\varphi^{-1}\bigg)
\nn\\
&=&{2}E_{23}(\varphi).\nn\label{II21non}
\end{eqnarray}
Here we have used the invariant (pseudo) inner product on
$gl(n,\mathbb R)$. Hence,
\begin{eqnarray}
\nabla E_{21}&=& \varphi_{xxxx}+\hat {\bf k}_1\label{II14non}
\end{eqnarray}
for some ${\bf k}$-valued smooth function $\hat {\bf k}_1$,
\begin{eqnarray}
\nabla E_{22}=
-\left(\varphi^{-1}\varphi_x\varphi^{-1}\varphi_x\varphi^{-1}\varphi_x\varphi^{-1}\right)_x+\hat{\bf
k}_2\label{II15non}
\end{eqnarray}
for some  ${\bf k}$-valued smooth function $\hat {\bf k}_2$ and
\begin{eqnarray}
\nabla
{E}_{23}=\frac{1}{2}\bigg(-(\varphi^{-1}\varphi_x\varphi^{-1}\varphi_x\varphi^{-1}\varphi_x\varphi^{-1})_x+
\hat{\bf k}_3\bigg)\label{II230non}
\end{eqnarray}
for some  ${\bf k}$-valued smooth function $\hat {\bf k}_3$. By
definition and the formulae (\ref{II14non},\ref{II15non}) and
(\ref{II230non}), the equation of the generalized bi-Schr\"odinger
flow of maps from ${\mathbb R}$ to ${\widetilde G}_{n,k}$ is
\begin{eqnarray}
\varphi_t&=&\bigg[\varphi,\,\alpha\nabla E+\beta\nabla E_2+\gamma\nabla {\widetilde E}\bigg]\nn\\
&=&\bigg[\varphi,\,\alpha\nabla E+\beta\nabla
\left(E_{21}-E_{22}\right)+
\left(\beta+2\gamma\right) \nabla E_{23}\bigg]\nn\\
&=&\bigg[\varphi,\,-{\alpha}\left(\varphi_{xx}-\varphi_x\varphi^{-1}\varphi_x\right)
+\beta\left(\varphi_{xxxx}+\hat{\bf
k}_1+(\varphi^{-1}\varphi_x\varphi^{-1}\varphi_x\varphi^{-1}\varphi_x\varphi^{-1})_x
+\hat{\bf k}_2\right)\nn\\
&&~~~
+\frac{1}{2}(\beta+2\gamma)\left(-(\varphi^{-1}\varphi_x\varphi^{-1}\varphi_x\varphi^{-1}\varphi_x\varphi^{-1})_x
+\hat{\bf
k}_3\right)\bigg]\nn\\
&=&\bigg[\varphi,\,-\alpha\varphi_{xx}
+\beta\varphi_{xxxx}+\left(-\gamma+\frac{\beta}{2}\right)
\varphi^{-1}\left(\varphi_x\varphi^{-1}\varphi_x\varphi^{-1}\varphi_x\right)_x\varphi^{-1}\bigg]
\nn\\
&=&\bigg[\varphi,\,-\alpha\varphi_{xx}
+\beta\varphi_{xxxx}+(4\gamma-2\beta)\left(\varphi_x\varphi^{-1}\varphi_x\varphi^{-1}\varphi_x\right)_x\bigg],\nn
\end{eqnarray}
which is still written as Eq.(\ref{bi0}). Here we have used the
identity (\ref{II16p}) in computation of the last equality. The
proof of Theorem \ref{Thm1} is completed. ~~$\Box$

\begin{Corollary}  \label{coro0}
For ${\bf g}=su(2)$, or equivalently $M=SU(2)/SU(1)\times
SU(1)=\{E^{-1}\sigma_3E~|~E\in SU(2)\}\cong\mathbb S^2$,
Eq.(\ref{bi0}) in Theorem \ref{Thm1} reduces exactly to the equation
(\ref{S2}) of the generalized bi-Schr\"{o}dinger flows from $\mathbb
R$ to $\mathbb S^2\hookrightarrow {\mathbb R}^{3}$.
\end{Corollary}
{\bf Proof}: In fact, for $E=\left(\begin{array}{cc}\bar a&-b\\\bar
b& a\end{array}\right)\in SU(2)$ (i.e. $|a|^2+|b|^2=1$), one sees
that
$$\varphi=\left(\begin{array}{cc}a&b\\-\bar b&\bar
a\end{array}\right)
\left(\begin{array}{cc}\frac{i}{2}&0\\0&-\frac{i}{2}\end{array}\right)
\left(\begin{array}{cc}\bar a&-b\\\bar b& a\end{array}\right)=
\frac{1}{2}\left(\begin{array}{cc}ix_1&x_2+ix_3\\-x_2+ix_3&
-ix_1\end{array}\right)
$$ with $x_1=|a|^2-|b|^2$, $x_2=-i(ab-\bar a\bar b)$, $x_3=-(ab+\bar a\bar b)$ and hence
$x_1^2+x_2^2+x_3^2=1$. Therefore, the correspondence between
$$M=SU(2)/SU(1)\times SU(1)=\{E^{-1}\sigma_3E~|~E\in SU(2)\}$$ and
$\mathbb S^2\hookrightarrow \mathbb R^3$ reads
$$
\varphi\longmapsto
\left(\begin{array}{c}x_1\\x_2\\x_3\end{array}\right)\in \mathbb
S^{2}\hookrightarrow\mathbb R^3.$$ In this case, one may verify
directly by using (\ref{II17}) and (\ref{II18}) that
\begin{eqnarray}
\varphi_x\varphi^{-1}\varphi_x\varphi^{-1} =
\left(\begin{array}{cc}-({x_1}_x)^2-({x_2}_x)^2-({x_3}_x)^2&0\\
0&-({x_1}_x)^2-({x_2}_x)^2-({x_3}_x)^2\end{array}\right).\nn\end{eqnarray}
Substituting this expression into Eq.(\ref{bi0}), one sees easily
that Eq.(\ref{bi0}) in this case is exactly the equation (\ref{S2})
of the generalized bi-Schr\"{o}dinger flows from $\mathbb R$ to the
2-sphere $\mathbb S^2\hookrightarrow \mathbb R^3$ given in \S2 (also
refer to \cite{DWW} or \cite{DLW}). ~~$\Box$

\begin{Corollary} \label{coro1}
For ${\bf g}=su(1,1)$, or equivalently $M=SU(1,1)/SU(1)\times
SU(1)=\{E^{-1}\sigma_3E~|~E\in SU(1,1)\}\cong\mathbb
H^2\hookrightarrow\mathbb R^{2,1}$, Eq.(\ref{bi0}) reduces exactly
to the equation (\ref{H2}) of the generalized bi-Schr\"{o}dinger
flows from $\mathbb R$ to $\mathbb H^2\hookrightarrow {\mathbb
R}^{2,1}$.
\end{Corollary}
{\bf Proof}: For $E=\left(\begin{array}{cc}\bar a&-b\\-\bar b&
a\end{array}\right)\in SU(2)$ (i.e. $|a|^2-|b|^2=1$), one sees that
$$\varphi=\left(\begin{array}{cc}a&b\\\bar b&\bar
a\end{array}\right)
\left(\begin{array}{cc}\frac{i}{2}&0\\0&-\frac{i}{2}\end{array}\right)
\left(\begin{array}{cc}\bar a&-b\\-\bar b& a\end{array}\right)=
\frac{1}{2}\left(\begin{array}{cc}ix_3&x_1+ix_2\\x_1-ix_1&
-ix_3\end{array}\right)
$$ with $x_3=|a|^2+|b|^2,\, x_2=-(ab+\bar a\bar b), \, x_1=-i(ab-\bar a\bar b)$ and hence
$x_1^2+x_2^2-x_3^2=-1$. Therefore, the correspondence between
$M=SU(1,1)/SU(1)\times SU(1)=\{E^{-1}\sigma_3E~|~E\in SU(2)\}$ and
$\mathbb H^2\hookrightarrow \mathbb R^{2,1}$ reads
$$
\varphi\longmapsto
\left(\begin{array}{c}x_1\\x_2\\x_3\end{array}\right)\in \mathbb
H^{2}\hookrightarrow\mathbb R^{2,1}.$$ In this case, one may verify
directly by using (\ref{II17}) and (\ref{II18}) that
\begin{eqnarray}
\varphi_x\varphi^{-1}\varphi_x\varphi^{-1} =
\left(\begin{array}{cc}({x_1}_x)^2+({x_2}_x)^2-({x_3}_x)^2&0\\
0&({x_1}_x)^2+({x_2}_x)^2-({x_3}_x)^2\end{array}\right).\nn\end{eqnarray}
Substituting this expression into Eq.(\ref{bi0}), one sees easily
that there is a sign difference in $\gamma$. This is because of our
new definition of the generalized bi-energy functional. By noting
this point, we see that Eq.(\ref{bi0}) in this case is exactly the
equation (\ref{H2}) of the generalized bi-Schr\"{o}dinger flows from
$\mathbb R$ to the hyperbolic plane $\mathbb H^2\hookrightarrow
\mathbb R^{2,1}$ given in \S2 (refer to \cite{DLW}). ~~$\Box$

\begin{Corollary}  \label{coro2}
For ${\bf g}=sl(2,\mathbb R)$, or equivalently $M=SL(2,\mathbb
R)/SL(1,\mathbb R)\times SL(1,\mathbb R)\cong\mathbb S^{1,1}$,
Eq.(\ref{bi0}) reduces exactly to the equation (\ref{S11}) of the
generalized bi-Schr\"{o}dinger flows from $\mathbb R$ to the de
Sitter 2-space $\mathbb S^{1,1}\hookrightarrow {\mathbb R}^{2,1}$.
\end{Corollary}
{\bf Proof}: Since a matrix $E\in SL(2,\mathbb R)$ is explicitly
expressed by  $$E=\left(\begin{array}{cc} a&-b\\-c&
d\end{array}\right)$$ with $ad-bc=1$, one sees that
$$\varphi=\left(\begin{array}{cc}a&b\\c&d\end{array}\right)
\left(\begin{array}{cc}\frac{1}{2}&0\\0&-\frac{1}{2}\end{array}\right)
\left(\begin{array}{cc} d&-b\\-c& a\end{array}\right)=
\frac{1}{2}\left(\begin{array}{cc}x_1&x_2+x_3\\x_2-x_3&
-x_1\end{array}\right)
$$ with $x_1=ad+bc, x_2=-(ab-cd), x_1=-(ab+cd)$ and hence
$x_1^2+x_2^2-x_3^2=1$. Therefore, the correspondence between
$M=SL(2,\mathbb R)/SL(1,\mathbb R)\times SL(1,\mathbb
R)=\{E^{-1}\sigma_3E~|~E\in SL(2,\mathbb R)\}$ and $\mathbb
S^{1,1}\hookrightarrow \mathbb R^{2,1}$ reads
$$
\varphi\longmapsto
\left(\begin{array}{c}x_1\\x_2\\x_3\end{array}\right)\in \mathbb
S^{1,1}\hookrightarrow \mathbb R^{2,1}.$$ In this case, one may
verify directly by using (\ref{II17}) and (\ref{II18}) that
\begin{eqnarray}
\varphi_x\varphi^{-1}\varphi_x\varphi^{-1} =
\left(\begin{array}{cc}-({x_1}_x)^2-({x_2}_x)^2+({x_3}_x)^2&0\\
0&-({x_1}_x)^2-({x_2}_x)^2+({x_3}_x)^2\end{array}\right).\nn\end{eqnarray}
Substituting this expression into Eq.(\ref{bi0}), one sees easily
that Eq.(\ref{bi0}) in this case is just the equation (\ref{S11}) of
the generalized bi-Schr\"{o}dinger flows from $\mathbb R$ to the de
Sitter 2-space $\mathbb S^{1,1}\hookrightarrow \mathbb R^{2,1}$
given in \S2 (also see \cite{DLW}). ~~$\Box$

\medskip
From Corollaries \ref{coro0}, \ref{coro1} and \ref{coro2} we see
that Eq.(\ref{bi0}) on the Hermitian symmetric Lie algebra $u(n)$
with index $k$ of compact type, and the Hermitian symmetric Lie
algebra $u(k,n-k)$ with index $k$ of noncompact type, and the
para-Hermitian symmetric Lie algebra $gl(n,\mathbb R)$ with index
$k$, is exactly the desired non-integrable third-order models of the
vortex filament on its corresponding symmetric Lie algebra. This is
because when $n=2$ and $k=1$, these models are equivalent
respectively to the Fukumoto-Moffatt's model (\ref{I0}) in $\mathbb
R^3$, the timelike third-order corrected model (\ref{I1}) and
spacelike third-order corrected model (\ref{I2}) of the vortex
filament in the Minkowski 3-space $\mathbb R^{2,1}$.

In the next section, like the Fukumoto-Moffatt's model (\ref{I0}) in
$\mathbb R^3$ and models (\ref{I1}) and (\ref{I2}) in $\mathbb
R^{2,1}$, we shall pursue whether Eq.(\ref{bi0}) on the Hermitian or
para-Hermitian symmetric Lie algebras admits the property expressed
as a motion of moving curves in the corresponding Lie algebra, and
whether Eq.(\ref{bi0}) on the Hermitian or para-Hermitian symmetric
Lie algebras can be transformed equivalently into a second to fourth
order nonlinear Schr\"odinger-like matrix equation.

\section{Sym-Pohlmeyer curves moving in symmetric Lie algebras}
In this section, we shall apply the technique of Sym (\cite{Sym})
and Pohlmeyer (\cite{Poh}) to produce geometric realizations of
Eq.(\ref{bi0}) in its corresponding symmetric Lie algebra ${\bf g}$.
Then, by applying the geometric concept of PDEs with the given
curvature representation proposed in the category of Yang-Mills
theory by the first author and his collaborators in \cite{ding1,
ding2} with the aid of gauge transformations, we transform
Eq.(\ref{bi0}) equivalently into a second to fourth order nonlinear
Schr\"odinger-like matrix equation. All these results return to what
are known for the Fukumoto-Moffatt's models (\ref{I0}), (\ref{I1})
and (\ref{I2}) of the vortex filament when ${\bf g}$ goes back to
$su(2)$, $su(1,1)$ and $sl(2,\mathbb R)$ respectively.

It should be mentioned that curves evolving in a symmetric Lie
algebra ${\bf g}$ described in \cite{Poh,Sym,LP} constitute a very
special subclass of the unit speed curves in ${\bf g}\cong {\mathbb
R}^{\hbox{dim}({\bf g})}$ (one notes here that ${\mathbb
R}^{\hbox{dim}({\bf g})}$ is regarded as the Euclidean space in the
case that ${\bf g}$ is in the first class or the pseudo-Euclidean
space otherwise). It is shown in \cite{LP} that the moving equation
(\ref{GLIE}): $\gamma_t=[\gamma_x,\gamma_{xx}]$ (i.e. the
leading-order model) in the symmetric Lie algebra $u(n)$ is a
geometric realization of the nonlinear Schr\"odinger matrix
equation. Meanwhile, the moving equation (\ref{GKDV}):
$\gamma_t=\gamma_{xxx}+\frac{3}{2}[\gamma_{xx},
[\gamma_x,\gamma_{xx}]]$ (i.e. the second-order model) in the
symmetric Lie algebra $sl(2n,\mathbb R)$ is a geometric realization
of the KdV type matrix equation (see \cite{dinghe}). Now the
following theorem shows that the third-order Eq.(\ref{bi0}) also
admits a geometric realization in its corresponding symmetric Lie
algebra.

\begin{Theorem}\label{Thm2}
The following moving equation in the Hermitian (or para-Hermitian)
symmetric Lie algebra ${\bf g}=$~$u(n)$ with index $k$, or
$u(k,n-k)$ with index $k$ or $gl(n,\mathbb R)$ with index $k$,
\begin{eqnarray}
{\widetilde\gamma}_t=-\alpha[{\widetilde\gamma}_x,\,{\widetilde\gamma}_{xx}]
+\beta\bigg([{\widetilde\gamma}_x,\,{\widetilde\gamma}_{xxxx}]
-[{\widetilde\gamma}_{xx},\,{\widetilde\gamma}_{xxx}]\bigg)+(4\gamma-2\beta)
\left[{\widetilde\gamma}_{x},\,
{\widetilde\gamma}_{xx}{\widetilde\gamma}_x^{-1}{\widetilde\gamma}_{xx}
{\widetilde\gamma}_x^{-1}{\widetilde\gamma}_{xx}\right]\label{GFM}
\end{eqnarray}
is equivalent to Eq.(\ref{bi0}) given in Theorem \ref{Thm1}, where
${\widetilde\gamma}\in {\bf g}$ with ${\widetilde\gamma}_x\in
\{E^{-1}\sigma_3E|E\in {\bf g}\}$.
\end{Theorem}

\begin{Remark} We point out that there are lots of moving Sym-Pohlmeyer curves
in ${\bf g}$ meeting the requirement of (\ref{GFM}). In fact, for a
solution $\varphi(t,x)=E^{-1}(t,x)\sigma_3E(t,x)$ to Eq.(\ref{bi0})
on ${\bf g}={\bf k}\oplus{\bf m}$ with the fixed index, letting
\begin{eqnarray}
\widetilde\gamma=\int^x_0\varphi(t,s)ds=\int^x_0E^{-1}(t,s)
\sigma_3E(t,s)ds, \nn\label{tilde}
\end{eqnarray}
we see that ${\widetilde\gamma}$ solves Eq.(\ref{GFM}) and satisfies
$({\widetilde\gamma})_x=\varphi(t,x)\in \{E^{-1}\sigma_3E|E\in {\bf
g}\}$. Theorem \ref{Thm2} indicates that the moving equation
(\ref{GFM}) is exactly the counterpart of the Fukumoto-Moffatt's
model on the symmetric Lie algebra ${\bf g}$, which is what we
pursue to find in this article.
\end{Remark}

{\bf Proof of Theorem \ref{Thm2}}: It is direct to verify, by taking
the derivative with respect to $x$, that if ${\widetilde\gamma}$
satisfies Eq.(\ref{GFM}) with ${\widetilde\gamma}_x\in
\{E^{-1}\sigma_3E|E\in {\bf g}\}$, then
$\varphi={\widetilde\gamma}_x$ satisfies Eq.(\ref{bi0}). Conversely,
if $\varphi(t,x)=E^{-1}(t,x)\sigma_3E(t,x)$ satisfies
Eq.(\ref{bi0}), then it is easy to see that
${\widetilde\gamma}=\int^x_0\varphi(t,s)ds$ solves Eq.(\ref{GFM})
and satisfies $({\widetilde\gamma})_x=\varphi(t,x)\in
\{E^{-1}\sigma_3E|E\in {\bf g}\}$. ~~$\Box$

\bigskip

It was proved by Hasimoto transform that the Fukumoto-Moffatt's
model (\ref{I0}) (resp. (\ref{I1},\ref{I2})) is equivalent to a
second to fourth order nonlinear Schr\"odinger-like equation (refer
to \cite{FMo,DLW}). In what follows, we shall show that
Eq.(\ref{GFM}) (or equivalently Eq.(\ref{bi0}))  is also equivalent
to a second to fourth order nonlinear Schr\"odinger-like matrix
equation. However, since we work on symmetric Lie algebras, the
Hasimoto transform is not applicable here. It is widely recognized
in literature that gauge transformation for a $1+1$ integrable
equation with zero curvature representation is equivalent to the
Hasimoto transform. Therefore, the gauge transformation is a good
candidate in higher dimensions that generalizes the Hasimoto
transform. But in our process we need to throw away the
integrability of the equations we deal with. Thanks to the works by
the first author and his collaborators in \cite{ding1, ding2} where
contrary to the concept of PDEs with zero curvature representation
(i.e. integrable PDEs), the geometric concept of PDEs with the given
curvature was proposed in the category of Yang-Mills theory. Thus,
we may still apply the gauge transformation to equations with a
given curvature representation that generalizes the Hasimoto
transform for these non-integrable PDEs.

Recall that by using the language in Yang-Mills theory,  the
generalized LIE (\ref{GLIE}) on the Hermitian symmetric Lie algebra
$u(n)$ with the index $k$ can be expressed by an equation with zero
curvature representation. More precisely, we first rewrite the
integrable Eq.(\ref{GLIE}) as
\begin{eqnarray}\varphi_t=[\varphi,\,\varphi_{xx}],\label{GDa}
\end{eqnarray}
where $\varphi={\widetilde \gamma}_x=E^{-1}\sigma_3E$ for some $E\in
U(n)$. Then we define a connection on the trivial bundle $\mathbb
R^2\times U(n)$ by
$$
A=\lambda\varphi dx+\left(\lambda^2\varphi+\lambda
[\varphi,\varphi_x]\right)dt.
$$
One may directly verify that $F_A=dA+A\wedge
A=\lambda\left(-\varphi_t+[\varphi,\varphi_{xx}]\right)$ and hence
Eq.(\ref{GDa}) is equivalent to the following zero curvature
representation:
$$
F_A=dA+A\wedge A=0.$$ For a given differential equation, if there
exists a connection 1-form $A$ and a given curvature 2-form $K$ such
that the equation is equivalent to the following formula:
$$F_A=dA+A\wedge A=K,$$
then it is called an equation with the given curvature
representation. We should point out that, roughly speaking,
equations with the given non-zero curvature representation are
non-integrable in general.

By using the above idea of the given curvature representation with
the aid of gauge transformations,  we may obtain the following
result.
\begin{Theorem} \label{Thm3}
On the symmetric Lie algebra ${\bf g}={\bf k}\oplus {\bf m}=$~$u(n)$
with index $k$, $u(k,n-k)$ with index $k$ or $gl(n,\mathbb R)$ with
index $k$, Eq.(\ref{bi0}) is equivalent to the following second to
fourth order nonlinear Schr\"odinger-like matrix equation by the
gauge transformation:
\begin{eqnarray}
2P_t\sigma_3&=&\alpha\bigg\{
-P_{xx}+2P^3\bigg\}+\beta\bigg\{P_{xxxx}-[P,[P,P_x]]_x-P\left((P^2)_{xx}-3P_x^3\right)\nn\\
&&-\left((P^2)_{xx}-3P_x^3\right)P-2(P^3)_{xx}+6P^5\bigg\}-2(8\gamma+\beta)
\bigg\{-(P^3)_{xx}+2P^5\nn\\
&&+P\left(\int^x_0(PPP_sP+PP_sPP)ds\right)
+\left(\int^x_0(PPP_sP+PP_sPP)ds\right)P\bigg\},\label{gauge}
\end{eqnarray}
where $P\in {\bf m}$ is an unknown ${\bf m}$-valued-matrix function.
\end{Theorem}
\noindent {\bf Proof}: 1) The case ${\bf g}=$~$u(n)$ with index $k$
($1\le k\le n-1$). First of all, we rewrite Eq.(\ref{bi0}) as
follows:
\begin{eqnarray}
\varphi_t =\bigg[\varphi,\,-\alpha\varphi_{xx}
+\beta\varphi_{xxxx}+4(4\gamma-2\beta)\left(\varphi^3_x\right)_x\bigg],\label{bi00}
\end{eqnarray}
since $\varphi_x\varphi^{-1}\varphi_x \varphi^{-1}=4\varphi_x^2$
(see (\ref{II20})). Next we come to construct a connection 1-form
$A$ and a curvature 2-form $K$ on the trivial bundle $\mathbb
R^2\times U(n)$ such that Eq.(\ref{bi00}) can be expressed by the
given curvature representation
\begin{eqnarray}
F_A=dA+A\wedge A=K. \label{bi14}
\end{eqnarray}
In fact, for a
given solution $\varphi=E^{-1}\sigma_3E$ to Eq.(\ref{bi00}) on ${\bf
g}=u(n)$ with index $k$, where the matrix $E\in U(n)$ is assumed to
satisfy $E_x=PE$ with $P\in {\bf m}$ as indicated in the proof of
Theorem \ref{Thm1}, we define
\begin{eqnarray}
A&=&\lambda \varphi dx+\bigg\{\lambda^4P_0+\lambda^3P_1+\lambda^2
P_2\nn\\
&&~+\lambda\left([\varphi,-\alpha\varphi_x+\beta\varphi_{xxx}+4(4\gamma-2\beta)\varphi^3_x]-\beta
[\varphi_x,\varphi_{xx}]\right)\bigg\}dt,\label{bi90}
\end{eqnarray}
where $P_0,P_1$ and $P_2$ independent of $\lambda$ are functions of
$P$ and its derivatives which will be determined later. Here the
choice of the coefficient of $\lambda$ in (\ref{bi90}) is to
guarantee the validity of Eq.(\ref{bi00}), which is enlightened by
the definition of the connection for the generalized LIE (\ref{GDa})
described above. It can be directly verified that
\begin{eqnarray}
F_A&=&dA+A\wedge
A\nn\\
&=&\bigg\{\lambda\bigg(-\varphi_t+\left[\varphi,-\alpha\varphi_x+\beta\varphi_{xxx}
+4(4\gamma-2\beta)\varphi^3_x\right]_x-\beta
[\varphi_x,\varphi_{xx}]_x\bigg)\nn\\
&&+\lambda^2\bigg(\left(P_2\right)_x+
\bigg[\varphi,\left[\varphi,-\alpha\varphi_x+\beta\varphi_{xxx}+4(4\gamma-2\beta)\varphi^3_x\right]-\beta
[\varphi_x,\varphi_{xx}]\bigg]
\bigg)\nn\\
&&+\lambda^3\bigg((\left(P_1\right)_x+ \left[\varphi,P_2\right]
\bigg)+\lambda^4\bigg(\left(P_0\right)_x+[\varphi,P_1]\bigg)\bigg\}dx\wedge
dt.\label{bi10}
\end{eqnarray}
By using the identities: $\varphi\varphi_{xxx}=-\varphi_{xxx}\varphi
-3(\varphi_x^2)_x$ and hence
\begin{eqnarray}
&&\left[\varphi,\left[\varphi,-\alpha\varphi_x+\beta\varphi_{xxx}+4(4\gamma-2\beta)\varphi^3_x\right]-\beta
[\varphi_x,\varphi_{xx}]\right]\nn\\
~~~~~~&&=~\alpha\varphi_x-\beta\varphi_{xxx}+6\beta(\varphi^2_x\varphi)_x-2(8\gamma+\beta)\varphi_x^3,\nn
\end{eqnarray}
which is deduced from the identities (\ref{II16}-\ref{II19}) and the
relation $E_x=PE$ with $P\in{\bf m}$, we see from (\ref{bi10}) that
in order for the two coefficients of $\lambda^2$ in both sides of
(\ref{bi14}) to be identically equal,
$$P_2=-\alpha\varphi+\beta(\varphi_{xx}-6\varphi^2_x\varphi)$$
and
\begin{eqnarray}
K=\bigg(-\lambda^22(8\gamma+\beta)\varphi_x^3\bigg)dx\wedge dt
\label{bi13}
\end{eqnarray}
is chosen to be the curvature 2-form. Furthermore, by the vanish of
the coefficients of $\lambda^3$ and $\lambda^4$ in the
left-hand-side of (\ref{bi14}) we obtain
$$
P_1=-\beta[\varphi,\varphi_x],~~P_0=-\beta\varphi$$ from
(\ref{bi10}). To summarize, if one defines
\begin{eqnarray}
A&=&\lambda \varphi
dx+\bigg\{-\lambda^4\beta\varphi-\lambda^3\beta[\varphi,\varphi_x]+\lambda^2
\left(-\alpha\varphi+\beta(\varphi_{xx}-6\varphi^2_x\varphi)\right)\nn\\
&&~+\lambda\left([\varphi,-\alpha\varphi_x+\beta\varphi_{xxx}+4(4\gamma-2\beta)\varphi^3_x]-\beta
[\varphi_x,\varphi_{xx}]\right)\bigg\}dt,\label{bi9}
\end{eqnarray}
and $K$ is given by (\ref{bi13}), then Eq.(\ref{bi00}) (or
equivalently Eq.(\ref{bi0})) is equivalent to (\ref{bi14}). In other
words, Eq.(\ref{bi00}) in this case possesses the given curvature
representation Eq.(\ref{bi14}) with $A$ and $K$ given by (\ref{bi9})
and (\ref{bi13}) respectively.

Now we shall make a gauge transformation to transform
Eq.(\ref{bi00}) into its equivalent second to fourth order nonlinear
Schr\"odinger-like matrix equation. Indeed, since the given solution
$\varphi=E^{-1}\sigma_3E$ to Eq.(\ref{bi00}) on ${\bf g}=u(n)$ with
index $k$  satisfies: $E_x=PE$ with $P\in {\bf m}$ as indicated
above, we make the following gauge transformation via $G=E^{-1}$:
\begin{eqnarray}
A\longmapsto {\widetilde A}=G^{-1}dG+G^{-1}AG \label{gauge1}
\end{eqnarray}
on the connection $A$ given by (\ref{bi9}). It is well-known from
the Yang-Mills theory that under the gauge transformation, the
relation between the curvature $F_A$ of $A$ and the curvature
$F_{\widetilde A}$ of $\widetilde A$ is:
\begin{eqnarray}
F_{\widetilde A}=G^{-1}F_AG=G^{-1}KG.\label{gauge2}
\end{eqnarray}

By a direct calculation, we see from (\ref{gauge1}) that
\begin{eqnarray}
{\widetilde A}&=&Ed(E^{-1})+EAE^{-1}=-(dE)E^{-1}+EAE^{-1}\nn\\
&=&\left(\lambda
\sigma_3-P\right)dx+\bigg\{-\lambda^4\beta\sigma_3+\lambda^3\beta
P+\lambda^2\left(-\alpha\sigma_3-2\beta(P_x+P^2)\sigma_3\right)\nn\\
&&+\lambda \left(\alpha
P+\beta(-P_{xx}-8P^3-[P,P_x])-4(4\gamma-2\beta)P^3\right)-E_t\bigg\}dt,\label{gauge3}
\end{eqnarray}
where $E_t$ independent of $\lambda$ is to be determined later. Here
we have used the identities (\ref{II16}-\ref{II20}). Furthermore, by
using (\ref{gauge3}) a direct computation shows that
\begin{eqnarray}
F_{\widetilde A}&=&d{\widetilde A}+{\widetilde A}\wedge{\widetilde
A}\nn\\
&=&\bigg\{P_t+\bigg(-\lambda^4\beta\sigma_3+\lambda^3\beta
P+\lambda^2\left(-\alpha\sigma_3-2\beta(P_x+P^2)\sigma_3\right)\nn\\
~~~~~&&~~+\lambda \left(\alpha
P+\beta(-P_{xx}-8P^3-[P,P_x])-4(4\gamma-2\beta)P^3\right)-E_t\bigg)_x\nn\\
~~~~~&&~~+\bigg[\lambda\sigma_3-P,-\lambda^4\beta\sigma_3+\lambda^3\beta
P+\lambda^2\left(-\alpha\sigma_3-2\beta(P_x+P^2)\sigma_3\right)\nn\\
~~~~~&&~~+\lambda \left(\alpha
P+\beta(-P_{xx}-8P^3-[P,P_x])-4(4\gamma-2\beta)P^3\right)-E_t\bigg]\bigg\}dx\wedge
dt\label{gauge4}
\end{eqnarray}
and
\begin{eqnarray}
{\widetilde
K}=G^{-1}F_AG=E\left(-\lambda^22(8\gamma+\beta)\varphi_x^3\right)E^{-1}dx\wedge
dt=\bigg(4(8\gamma+\beta)P^3\sigma_3\bigg)dx\wedge dt.\label{gauge5}
\end{eqnarray}
Substituting (\ref{gauge4}) and (\ref{gauge5}) into (\ref{gauge2}),
we see that the coefficients of $\lambda^5$, $\lambda^4$,
$\lambda^3$ and $\lambda^2$ in both sides of (\ref{gauge2}) are
equal to each other. The fact that the  coefficient of $\lambda$ in
the left-hand-side of (\ref{gauge2}) is zero implies that
$$[\sigma_3,E_t]=\alpha P_x-\beta
(P_{xxx}-[P,[P,P_x]])-16\gamma (P^3)_x$$ and hence
\begin{eqnarray}
E_t^{(\hbox{off-diag})}=2\alpha P_x\sigma_3-2\beta
(P_{xxx}-[P,[P,P_x]])\sigma_3-32\gamma
(P^3)_x\sigma_3,\label{gauge6}
\end{eqnarray}
where $E_t^{(\hbox{off-diag})}$ is the ${\bf m}$-part or
off-diagonal part of $E_t$. Here we have used the fact: for a matrix
$X\in u(n)={\bf k}\oplus{\bf m}$,
$$[\sigma_3,X]\sigma_3=\frac{1}{2}X^{(\hbox{off-diag})}.$$
Meanwhile, the fact that the coefficient of $\lambda^0$ in the
left-hand-side of (\ref{gauge2}) is zero implies that
\begin{eqnarray}
P_t-(E_t)_x+[P,E_t]=0.\label{gauge7}
\end{eqnarray}
By taking the diagonal part of (\ref{gauge7}), we see that
\begin{eqnarray}
\left(E_t^{(\hbox{diag})}\right)_x&=&\left[P,E_t^{(\hbox{off-diag})}\right]=\left[P,2\alpha
P_x\sigma_3-2\beta (P_{xxx}-[P,[P,P_x]])\sigma_3-32\gamma
(P^3)_x\sigma_3\right]\nn\\
&=&2\alpha
\left(P^2\right)_x\sigma_3-2\beta\left((P^2)_{xx}-3(P^2_x)\right)_x\sigma_3
+6\beta\left(P^4\right)_x\sigma_3\nn\\&&
-(8\gamma+\beta)\left(4(P^4)_x+4(PPP_xP+PP_xPP)\right)\sigma_3\nn
\end{eqnarray} and hence
\begin{eqnarray}
E_t^{(\hbox{diag})} &=&2\alpha
P^2\sigma_3-2\beta\left((P^2)_{xx}-3(P^2_x)\right)\sigma_3 +6\beta
P^4\sigma_3\nn\\&&
-(8\gamma+\beta)\left(4P^4+4\int^x_0(PPP_sP+PP_sPP)ds\right)\sigma_3.\label{gauge8}
\end{eqnarray}
By taking the off-diagonal part of (\ref{gauge7}):
$$P_t-\left(E_t^{(\hbox{off-diag})}\right)_x+\left[P,E_t^{(\hbox{diag})}\right]=0$$
and using (\ref{gauge6},\ref{gauge8}), we finally arrive at the
following second to fourth order nonlinear Schr\"odinger matrix
equation:
\begin{eqnarray}
2P_t\sigma_3&=&\alpha\bigg\{
-P_{xx}+2P^3\bigg\}+\beta\bigg\{P_{xxxx}-[P,[P,P_x]]_x-P\left((P^2)_{xx}-3P_x^3\right)\nn\\
&&-\left((P^2)_{xx}-3P_x^3\right)P-2(P^3)_{xx}+6P^5\bigg\}-2(8\gamma+\beta)
\bigg\{-(P^3)_{xx}+2P^5\nn\\
&&+P\left(\int^x_0(PPP_sP+PP_sPP)ds\right)
+\left(\int^x_0(PPP_sP+PP_sPP)ds\right)P\bigg\},\nn\label{gauge9}
\end{eqnarray}
which is just Eq.(\ref{gauge}). This proves that any solution to
Eq.(\ref{bi00}) is transformed into a solution to Eq.(\ref{gauge})
by the gauge transformation. Conversely, it is easy to see that
every step in the above proof is invertible and hence any solution
to Eq.(\ref{gauge}) is transformed into a solution to
Eq.(\ref{bi00}). The proof that Eq.(\ref{bi0}) on $u(n)$ is
equivalent to Eq.(\ref{gauge}) is completed.

\medskip

2) The case ${\bf g}=$~$u(k,n-k)$ with index $k$ ($1\le k\le n-1$).
Since the identities (\ref{II16}-\ref{II20}) are the same as in the
case ${\bf g}=$~$u(n)$ with index $k$, following exactly the line in
the proof of 1), we see that Eq.(\ref{bi0}) on $u(k,n-k)$ with the
index $k$ is transformed into Eq.(\ref{gauge}) and vice versa. This
proves that Eq.(\ref{bi0}) on $u(k,n-k)$ is equivalent to
Eq.(\ref{gauge}).

\medskip

3) The case ${\bf g}=$~$gl(n,\mathbb R)$ with index $k$ ($1\le k\le
n-1$). First of all, we rewrite Eq.(\ref{bi0}) in this case as
follows:
\begin{eqnarray}
\varphi_t =\bigg[\varphi,-\alpha\varphi_{xx}
+\beta\varphi_{xxxx}-4(4\gamma-2\beta)\left(\varphi^3_x\right)_x\bigg],\label{bi00p}
\end{eqnarray}
since $\varphi_x\varphi^{-1}\varphi_x \varphi^{-1}=-4\varphi_x^2$
(see (\ref{II20p})). Next we come to construct a connection 1-form
$A$ and a curvature 2-form $K$ on the trivial bundle $\mathbb
R^2\times gl(n,\mathbb R)$ such that Eq.(\ref{bi00p}) can be
expressed as an equation with the given curvature representation. In
fact, unlike in the proof of Theorem \ref{Thm1}, we replace a
solution of the form $\varphi=E^{-1}\sigma_3E$ to Eq.(\ref{bi00p})
by $\varphi=E\sigma_3E^{-1}$ with the matrix $E\in GL(n,\mathbb R)$
being now assumed to satisfy $E_x=EP$ with $P\in {\bf m}$. One notes
that the new $P$ here equals  $-P$ with $P$ used in the proof of
Theorem \ref{Thm1}. By using the relation $E_x=EP$ and hence the
identity:
\begin{eqnarray}
&&\left[\varphi,\left[\varphi,-\alpha\varphi_{xx}
+\beta\varphi_{xxxx}-4(4\gamma-2\beta)\left(\varphi^3_x\right)_x\right]\right]
-\beta\left[\varphi,[\varphi_x,\varphi_{xx}]\right]
\nn\\&=&-\alpha\varphi_x+\beta\left(\varphi_{xxx}+6(\varphi_{x}^2)_x\varphi\right)-
4(4\gamma-2\beta)\varphi^3_x-4\beta\varphi_x^3,\nn
\end{eqnarray}
we defines a connection 1-form by
\begin{eqnarray}
A&=&\lambda \varphi
dx+\bigg\{-\lambda^4\beta\varphi+\lambda^3\beta[\varphi,\varphi_x]+\lambda^2
\left(\alpha\varphi-\beta(\varphi_{xx}+6\varphi^2_x\varphi)\right)\nn\\
&&~+\lambda\left([\varphi,-\alpha\varphi_x+\beta\varphi_{xxx}-4(4\gamma-2\beta)\varphi^3_x]-\beta
[\varphi_x,\varphi_{xx}]\right)\bigg\}dt,\label{bi9p}
\end{eqnarray}
and
\begin{eqnarray}
K=\bigg(-\lambda^22(8\gamma+\beta)\varphi_x^3\bigg)dx\wedge dt
\label{bi13p}
\end{eqnarray}
to be a given curvature 2-form. Then Eq.(\ref{bi00p}) (or
equivalently Eq.(\ref{bi0}) on $gl(n,\mathbb R)$) admits the
following given curvature representation:
\begin{eqnarray}
F_A=dA+A\wedge A=K.\label{bi14p}
\end{eqnarray}

The key point here in making a gauge transformation for
Eq.(\ref{bi14p}) is as follows. Since the given solution
$\varphi=E\sigma_3E^{-1}$ to Eq.(\ref{bi00p}) on ${\bf
g}=gl(n,\mathbb R)$ with index $k$ is obtained by assuming that
$E\in GL(n,\mathbb R)$ satisfies $E_x=EP$ with $P\in {\bf m}$. The
gauge transformation via $G=E$ for the connection $A$ given by
(\ref{bi9p}) is:
\begin{eqnarray}
A\longmapsto {\widetilde A}=E^{-1}dE+E^{-1}AE. \label{gauge1p}
\end{eqnarray}
Hence
\begin{eqnarray}
{\widetilde K}=F_{\widetilde A}=E^{-1}F_AE=E^{-1}KE=\bigg(
\lambda^24(8\gamma+\beta)P^3\sigma_3\bigg)dx\wedge
dt.\nn\label{gauge2p}
\end{eqnarray}

It is also through a direct computation and by using
(\ref{II16p}-\ref{II20p}) that ${\widetilde A}$ given by
(\ref{gauge1p}) can be explicitly expressed as
\begin{eqnarray}
{\widetilde A }&=&\left(\lambda
\sigma_3+P\right)dx+\bigg\{-\lambda^4\beta\sigma_3-\lambda^3\beta
P+\lambda^2\left(\alpha\sigma_3-2\beta(P_x-P^2)\sigma_3\right)\nn\\
&&+\lambda \left(\alpha P-\beta(P_{xx}-[P,P_x])-16\gamma
P^3\right)-E_t\bigg\}dt,\label{gauge3p}
\end{eqnarray}
where $E_t$ independent of $\lambda$ is to be determined later.
Similarly, from the relation $F_{\widetilde A}=d{\widetilde A}+
{\widetilde A}\wedge{\widetilde A}={\widetilde K}$, one may obtain
$$E^{(\hbox{off-diag})}_t=-2\alpha P_x\sigma_3+2\beta
P_{xxx}\sigma_3-2\beta[P,[P,P_x]]\sigma_3+32\gamma (P^3)_x\sigma_3$$
and \begin{eqnarray} \left(E^{(\hbox{diag})}_t\right)_x&=&2\alpha
(P^2)_x\sigma_3-2\beta
\left((P^2)_{xx}-3P^2_x\right)_x\sigma_3+6\beta\left(P^4\right)_x\sigma_3\nn\\
&&-(32\gamma
+4\beta)\left((P^4)_x+P^2P_xP+PP_xP^2\right)\sigma_3.\nn
\end{eqnarray}
The above identities come from equalizing the coefficients of
$\lambda$ and the diagonal part of the constant terms in
$d{\widetilde A}+ {\widetilde A}\wedge{\widetilde A}={\widetilde K}$
respectively. Based on the above two identities, the off-diagonal
part of the constant term in $d{\widetilde A}+ {\widetilde
A}\wedge{\widetilde A}={\widetilde K}$ implies
$$-P_t+\left(E^{(\hbox{off-diag})}_t\right)_x+\left[P,E^{(\hbox{diag})}_t\right]=0,$$
which can be explicitly rewritten as
\begin{eqnarray}
2P_t\sigma_3&=&\alpha\bigg\{
-P_{xx}+2P^3\bigg\}+\beta\bigg\{P_{xxxx}-[P,[P,P_x]]_x-P\left((P^2)_{xx}-3P_x^3\right)\nn\\
&&-\left((P^2)_{xx}-3P_x^3\right)P-2(P^3)_{xx}+6P^5\bigg\}-2(8\gamma+\beta)
\bigg\{-(P^3)_{xx}+2P^5\nn\\
&&+P\left(\int^x_0(PPP_sP+PP_sPP)ds\right)
+\left(\int^x_0(PPP_sP+PP_sPP)ds\right)P\bigg\}. \nn\label{gauge9p}
\end{eqnarray}
This equation is exactly Eq.(\ref{gauge}) and hence any solution to
Eq.(\ref{bi00p}) is transformed to a solution to Eq.(\ref{gauge}) by
the gauge transformation. One also notes that all the steps in the
above proof are invertible, or in other words, Eq.(\ref{gauge}) in
this case is transformed to Eq.(\ref{bi00p}) by the gauge
transformation. The proof that Eq.(\ref{bi00p}) on $gl(n, \mathbb
R)$ and Eq.(\ref{gauge12}) are equivalent to each other is finished.
~~ $\Box$

\medskip
The following proposition gives a detailed description of
Eq.(\ref{gauge}) on the three different types of symmetric Lie
algebras.

\begin{Proposition}\label{Thm4}
Eq.(\ref{gauge}) given in Theorem \ref{Thm3} on the Hermitian
symmetric Lie algebra $u(n)$ with index $k$ of compact type, the
Hermitian symmetric Lie algebra $u(k,n-k)$ with index $k$ of
noncompact type and the para-Hermitian symmetric Lie algebra
$gl(n,\mathbb R)$ with index $k$ are respectively given by
\begin{eqnarray}
&&iq_t+\alpha\bigg\{-q_{xx}-2qq^*q\bigg\}+\beta\bigg\{q_{xxxx}+4q_{xx}q^*q+2qq^*_{xx}q+4qq^*q_{xx}
\nn\\
&&~~~~+2q_xq^*_xq+6q_xq^*q_x+2qq^*_xq_x+6qq^*qq^*q\bigg\}-2(8\gamma+\beta)\bigg\{(qq^*q)_{xx}+
2qq^*qq^*q\nn\\
&&~~~~+q\left(\int^xq^*(qq^*)_sqds\right)
+\left(\int^xq(q^*q)_sq^*ds\right)q\bigg\}=0, \label{gauge10}
\end{eqnarray}
where $q$ is a $k\times(n-k)$ complex-matrix-valued unknown function
and $q^*$ stands for the transposed conjugate matrix of $q$,
\begin{eqnarray}
&&iq_t+\alpha\bigg\{-q_{xx}+2qq^*q\bigg\}+\beta\bigg\{q_{xxxx}-4q_{xx}q^*q-2qq^*_{xx}q-4qq^*q_{xx}
\nn\\
&&~~~~-2q_xq^*_xq-6q_xq^*q_x-2qq^*_xq_x+6qq^*qq^*q\bigg\}-2(8\gamma+\beta)\bigg\{-(qq^*q)_{xx}+
2qq^*qq^*q\nn\\
&&~~~~+q\left(\int^xq^*(qq^*)_sqds\right)
+\left(\int^xq(q^*q)_sq^*ds\right)q\bigg\}=0, \label{gauge11}
\end{eqnarray}
where $q$ is a $k\times(n-k)$ complex-matrix-valued unknown
function, and
\begin{eqnarray}
\left\{\begin{aligned}
q_t&=\alpha\bigg\{q_{xx}-2qrq\bigg\}+\beta\bigg\{-q_{xxxx}
+4q_{xx}rq+2qr_{xx}q+4qrq_{xx}+2q_xr_xq\\
&~~~~+6q_xrq_x+2qr_xq_x
-6qrqrq\bigg\}+2(8\gamma+\beta)\bigg\{-(qrq)_{xx}+
2qrqrq\\
&~~~~+q\left(\int^xr(qr)_sqds\right)
+\left(\int^xq(rq)_srds\right)q\bigg\},\\
r_t&=\alpha\bigg\{-r_{xx}+2rqr\bigg\}+\beta\bigg\{r_{xxxx}
-4r_{xx}qr-2rrq_{xx}r-4rqr_{xx}-2r_xq_xr\\
&~~~~-6r_xqr_x-2rq_xr_x
+6rqrqr\bigg\}+2(8\gamma+\beta)\bigg\{(rqr)_{xx}-
2rqrqr\\
&~~~~-r\left(\int^xq(rq)_srds\right)
-\left(\int^xr(qr)_sqds\right)r\bigg\},\label{gauge12}
\end{aligned}\right.
\end{eqnarray}
where $q$ is a $k\times(n-k)$ real-matrix-valued unknown function
and $r$ an $(n-k)\times k$ real-matrix-valued unknown function.
\end{Proposition}
\noindent {\bf Proof}: 1) In the case ${\bf g}={\bf k}\oplus{\bf
m}=u(n)$ with index $k$, we write $P\in {\bf m}\subset u(n)$ as
\begin{eqnarray}
P=\left(\begin{array}{cc}0&q\\
-q^*&0\end{array}\right), \nn
\end{eqnarray}
where $q$ is a $k\times(n-k)$ complex matrix and $q^*$ stands for
the transposed conjugate matrix of $q$, and then by the identities
(\ref{II16}-\ref{II20}) and a long direct computation, we see that
Eq.(\ref{gauge9}) is equivalent to
\begin{eqnarray}
\left(\begin{array}{cc}0&-iq_t\\
-iq^*_t&0\end{array}\right)&=&\alpha \left(\begin{array}{cc}0&-q_{xx}-2qq^*q\\
q^*_{xx}+2q^*q^*&0\end{array}\right)+\beta\left(\begin{array}{cc}0&A\\
-A^*&0\end{array}\right)\nn\\
&&-2(8\gamma+\beta)\left(\begin{array}{cc}0&B\\
-B^*&0\end{array}\right),\nn
\end{eqnarray}
where
$$
A=q_{xxxx}+4q_{xx}q^*q+2qq^*_{xx}q+4qq^*q_{xx}+2q_xq^*_xq+6q_xq^*q_x+2qq^*_xq_x+6qq^*qq^*q
$$and
 $$B=(qq^*q)_{xx}+2qq^*qq^*q+q\left(\int^x_0q^*(qq^*)_sqds\right)
+\left(\int^xq^*(q^*q)_sq^*ds\right)q.$$ Therefore, Eq.(\ref{gauge})
in this case can be explicitly expressed by the following second to
fourth order nonlinear Schr\"odinger matrix equation:
\begin{eqnarray}
&&iq_t+\alpha\bigg\{-q_{xx}-2qq^*q\bigg\}+\beta\bigg\{q_{xxxx}+4q_{xx}q^*q+2qq^*_{xx}q+4qq^*q_{xx}
\nn\\
&&~~~~+2q_xq^*_xq+6q_xq^*q_x+2qq^*_xq_x+6qq^*qq^*q\bigg\}-2(\gamma+\beta)\bigg\{(qq^*q)_{xx}+
2qq^*qq^*q\nn\\
&&~~~~+q\left(\int^xq^*(qq^*)_sqds\right)
+\left(\int^xq^*(q^*q)_sq^*ds\right)q\bigg\}=0, \nn
\end{eqnarray}
which is just Eq.(\ref{gauge10}).

2) In the case ${\bf g}={\bf k}\oplus{\bf m}=u(k,n-k)$ with index
$k$, we now write $P\in{\bf m}\subset u(k,n-k)$ as
\begin{eqnarray}
P=\left(\begin{array}{cc}0&q\\
q^*&0\end{array}\right), \nn
\end{eqnarray}
where $q$ is a $k\times(n-k)$ complex matrix and $q^*$ stands for
the transposed conjugate matrix of $q$. Following the same
computation as above, we see that Eq.(\ref{gauge}) on $u(k,n-k)$ is
explicitly expressed by Eq.(\ref{gauge11}).

3) In the case ${\bf g}={\bf k}\oplus{\bf m}=gl(n,\mathbb R)$ with
index $k$, we now write $P\in {\bf m}\subset gl(n,\mathbb R)$ as
\begin{eqnarray}
P=\left(\begin{array}{cc}0&q\\
r&0\end{array}\right), \nn
\end{eqnarray}
where $q$ is a $k\times(n-k)$ unknown real-valued-matrix and $r$ is
an $(n-k)\times k$ unknown real-valued-matrix. By using the
identities (\ref{II16p}-\ref{II20p}) and a long computation, one may
verify directly that Eq.(\ref{gauge}) is explicitly expressed by
Eq.(\ref{gauge12}). The details are omitted here. The proof of the
Proposition is completed.~~$\Box$

\begin{Remark}\label{Remark2}
Unexpectedly, from  Theorem \ref{Thm3} or Proposition \ref{Thm4}
Eqs.(\ref{gauge10},\ref{gauge11}) and (\ref{gauge12}) are
differential-integral equations in general. A peculiar feature of
them is that these equations reduce to the second to fourth order
completely integrable nonlinear Schr\"odinger AKNS-matrix
differential equations when $\gamma=-\frac{\beta}{8}$, and
especially when $\alpha=0$ and $\gamma=-\frac{\beta}{8}$,
Eq.(\ref{gauge10}) goes back to the fourth order integrable matrix
Schr\"odinger equation (\ref{AKNS01}) on $u(n)$ given in \S2. Hence
the three typical second to fourth order integrable nonlinear
Schr\"odinger AKNS-matrix equations also have geometric realizations
of Sym-Pohlmeyer curves satisfying Eq.(\ref{GFM}) with
$\gamma=-\frac{\beta}{8}$ in the symmetric Lie algebra ${\bf g}={\bf
k}\oplus{\bf m}=u(n)$, $u(k,n-k)$ and $gl(n,\mathbb R)$
respectively.
\end{Remark}

Another distinguishing feature of Eqs.(\ref{gauge10},\ref{gauge11})
and (\ref{gauge12}) is that these differential-integral equations
reduce to partial differential equations when $n=2$.

\begin{Corollary} \label{Cor3}
On ${\bf g}=su(2)$, Eq.(\ref{gauge10}) reduces exactly to
\begin{eqnarray}
&&iq_t-\alpha\bigg\{q_{xx}+2|q|^2q\bigg\}+\beta\bigg\{
q_{xxxx}+6\bigg[|q|^2q_{xx}+(q_x)^2\bar q\bigg]\nonumber\\
~~~&&+\left[6|q|^4+2(|q|^2)_{xx}\right]
q\bigg\}-2\left(8\gamma+{\beta}\right)\bigg\{
(|q|^2q)_{xx}+3|q|^4q\bigg\}=0, \label{G10}
\end{eqnarray}
where $q$ is an unknown complex valued function. This equation is
just the second to fourth order equation that equivalent to the
Fukumoto-Maffott's model (\ref{I0}) or in other words the
generalized Schr\"{o}dinger flows from $\mathbb R$ to the 2-sphere
$\mathbb S^{2}\hookrightarrow R^3$ (refer to \cite{FMo,DWW1}).
\end{Corollary}
\noindent {\bf Proof}: In this case, $q^*=\bar q$ with $q$ being an
unknown complex function of $t$ and $x$. One may verify that
\begin{eqnarray}
&&q_{xxxx}+4q_{xx}q^*q+2qq^*_{xx}q+4qq^*q_{xx}
+2q_xq^*_xq+6q_xq^*q_x+2qq^*_xq_x+6qq^*qq^*q\nn\\&&~=~~
q_{xxxx}+6\bigg[|q|^2q_{xx}+(q_x)^2\bar
q\bigg]+\left[6|q|^4+2(|q|^2)_{xx}\right] q\nn
\end{eqnarray}
and (up to a term $c(t)q$ with $c(t)$ depending only on $t$ which
can be concealed by the transform of the form: $q\mapsto
q\exp\left(i\int^tc(s)ds\right)$)
$$q\left(\int^xq^*(qq^*)_sqds\right)
+\left(\int^xq^*(q^*q)_sq^*ds\right)q=|q|^4q.$$ Substituting the
above into (\ref{gauge10}), we see that Eq.(\ref{gauge10}) becomes
exactly Eq.(\ref{G10}).~~$\Box$

\medskip

Similarly, we have
\begin{Corollary}\label{Cor4}
On ${\bf g}=su(1,1)$, Eq.(\ref{gauge11}) reduces exactly to
\begin{eqnarray}
&&iq_t-\alpha\bigg\{q_{xx}-2|q|^2q\bigg\}+\beta\bigg\{
q_{xxxx}-6\bigg[|q|^2q_{xx}+(q_x)^2\bar q\bigg]\nonumber\\
~~~&&+\left[6|q|^4-2(|q|^2)_{xx}\right]
q\bigg\}-2\left(8\gamma+{\beta}\right)\bigg\{
(|q|^2q)_{xx}-3|q|^4q\bigg\}=0, \label{G11}
\end{eqnarray}
where $q$ is an unknown complex valued function. This equation is
just the second to fourth order equation that equivalent to the
timelike third-order corrected model (\ref{I1}) of the vortex
filament in $\mathbb R^{2,1}$ or in other words the generalized
generalized Schr\"{o}dinger flows from $\mathbb R$ to the hyperbolic
2-space $\mathbb H^{2}\hookrightarrow R^{2,1}$ (see
\cite{DWW1,DLW}).
\end{Corollary}

\begin{Corollary} \label{Cor5}
On ${\bf g}=sl(2,\mathbb R)$, Eq.(\ref{bi0}) reduces exactly to
\begin{eqnarray}
\left\{\begin{aligned} q_t&=\alpha (q_{xx}-2q^2r)+\beta\bigg(
-q_{xxxx}+6q^2_xr+4qq_xr_x+8qrq_{xx}+2q^2r_{xx}-6q^3r^2\bigg)
\\&~~~~~~~~~~~~~~~~~~-2(8\gamma+\beta)\bigg((qrq)_{xx}-3q^3r^2\bigg),\\
r_t&=-\alpha (r_{xx}-2qr^2)+\beta\bigg(
r_{xxxx}-6qr_x^2-4rq_xr_x-8qrr_{xx}-2r^2q_{xx}+6q^2r^3\bigg)
\\&~~~~~~~~~~~~~~~~~+2(8\gamma+\beta)\bigg((rqr)_{xx}-3q^2y^3\bigg),\label{G12}
\end{aligned}\right.
\end{eqnarray}
where $q$ and $r$ are unknown real functions. This equation is just
the second to fourth order equation that equivalent to the spacelike
third-order corrected model (\ref{I2}) of the vortex filament in
$\mathbb R^{2,1}$ or equivalently  the generalized Schr\"{o}dinger
flows from $\mathbb R$ to the de Sitter 2-space $\mathbb
S^{1,1}\hookrightarrow R^{2,1}$ (refer to \cite{DLW}).
\end{Corollary}

From Proposition \ref{Thm4} and Corollaries \ref{Cor3}, \ref{Cor4}
and \ref{Cor5}, we see that our exploitation of the third-order
models of the vortex filament on  Hermitian or para-Hermitian
symmetric Lie algebras ${\bf g}$ returns exactly to what we knew in
literature for the Fukumoto-Moffatt's models (\ref{I0}) in the
Euclidean 3-space $\mathbb R^3$ or the models (\ref{I1}) and
(\ref{I2}) of the vortex filament in the Minkowski 3-space $\mathbb
R^{2,1}$  when ${\bf g}=su(2)$, or $su(1,1)$, or $sl(2,\mathbb R)$.

\section{Conclusions and remarks}

In this article, by modifying the concept of generalized
bi-Schr\"odinger maps into K\"ahler or para-K\"ahler manifolds, we
have established the third-order non-integrable model of the vortex
filament on a Hermitian or para-Hermitian symmetric Lie algebra in a
unified way.  Combining it with what are known in literature for the
leading order and the second order models, we may list all the basic
models of the vortex filament on a symmetric Lie algebra ${\bf g}$
up to the third-order approximation, which are moving equations of
Sym-Pohlmeyer's curves in the symmetric Lie algebra ${\bf g}$:
\begin{eqnarray}
{\widetilde\gamma}_t&=&[{\widetilde\gamma}_x,\,{\widetilde\gamma}_{xx}], \label{GLIE1}\\
{\widetilde\gamma}_t&=&{\widetilde\gamma}_{xxx}+\frac{3}{2}[{\widetilde\gamma}_{xx},\,
[{\widetilde\gamma}_x,\,{\widetilde\gamma}_{xx}]],\label{GKDV1}\\
{\widetilde\gamma}_t&=&\beta\bigg([{\widetilde\gamma}_x,\,{\widetilde\gamma}_{xxxx}]
-[{\widetilde\gamma}_{xx},\,{\widetilde\gamma}_{xxx}]\bigg)+(4\gamma-2\beta)
\left[{\widetilde\gamma}_{x},\,
{\widetilde\gamma}_{xx}{\widetilde\gamma}_x^{-1}{\widetilde\gamma}_{xx}
{\widetilde\gamma}_x^{-1}{\widetilde\gamma}_{xx}\right],\label{GFM1}
\end{eqnarray}
where $[\cdot,\cdot]$ denotes the Lie bracket on ${\bf g}$.
Comparing these equations with the models existing in the localized
induction (matrix) hierarchy, we see that Eq.(\ref{GLIE1}) (resp.
Eq.(\ref{GKDV1})) coincides exactly with Eq.(\ref{GLIE}) (resp.
Eq.(\ref{GKDV})). However, Eq.(\ref{GFM1}) is not the fourth order
equation in the localized induction (matrix) hierarchy in general,
unless the two parameters $\beta$ and $\gamma$ satisfy the
condition: $\gamma=-\frac{\beta}{8}$ as mentioned in Remark
\ref{Remark2}. We should point out that the discovery of the
third-order correction model (\ref{GFM1}) is based on exploiting the
intrinsic geometry of the Hermitian or para-Hermitian symmetric Lie
algebras and the corresponding symmetric spaces, especially their
K\"ahler or para-K\"ahler structures. The application of generalized
bi-Schr\"odinger flows overcomes the limitation of the methods in
the theory of integrable systems and creates the model in a purely
geometric way. Furthermore, by using the concept of equations with a
given curvature representation with the aid of gauge
transformations, the models (\ref{GFM1}) of the vortex filament on
the three types of symmetric Lie algebras are transformed to three
types of fourth order nonlinear Schr\"odinger-like matrix
differential-integral equations. These three types of equations
reduce to the three typical 4th order integrable nonlinear
Schr\"odinger matrix differential equations in the AKNS matrix
hierarchy when the parameters $\beta$ and $\gamma$ satisfy the
relation: $\gamma=-\frac{\beta}{8}$. It is very appealing and
delighting that when the symmetric Lie algebra ${\bf g}$ goes back
to $su(2)$, $su(1,1)$ and $sl(2,\mathbb R)$ respectively, the models
(\ref{GLIE1},\ref{GKDV1}) and (\ref{GFM1}) on ${\bf g}$ reduce
respectively to the LIE-type models (\ref{Da},\ref{New041}) and
(\ref{New042}), the Fukumoto-Miyazaki-type models
(\ref{005},\ref{ab1}) and (\ref{ab2}), the Fukumoto-Moffatt-type
models (\ref{I0},\ref{I1}) and (\ref{I2}) in $\mathbb R^3$ or
$\mathbb R^{2,1}$, and all the theory we know about
Eqs.(\ref{GLIE1},\ref{GKDV1}) and (\ref{GFM1}) are reduced to the
corresponding theory of the vortex filament in $\mathbb R^3$ or
$\mathbb R^{2,1}$.

This indicates that there exists a complete hidden theory of the
vortex filament on symmetric Lie algebras, which consists of motions
of arclength-parameterized Sym-Pohlmeyer curves evolving in the
(Hermitian or para-Hermitian) symmetric Lie algebras we deal with.
Models (\ref{GLIE1},\ref{GKDV1}) and (\ref{GFM1}) provide the basic
equations of the vortex filament on the symmetric Lie
fluid algebra up to the third-order approximation. The theory of the
vortex filament in the Euclidean 3-space $\mathbb R^3$ or in the
Minkowski 3-space $\mathbb R^{2,1}$ can be understood uniformly on
the level of symmetric Lie algebras.

The study of the vortex filament in the Euclidean and Minkowski
3-space or on symmetric Lie algebras is always significant. We
believe that the model (\ref{GFM1}) and its related equations
(\ref{bi0}) and (\ref{gauge}) as well as the geometric concept of
generalized bi-Schr\"odinger flows will have much more deeper
applications in both mathematics and physics. However, there are
some questions still unclear at the present time. For example, 1)
The para-K\"ahler structure in geometry seems to be peculiar
beautiful and mysterious to us, which has not been well understood
up to now. How to characterize the dynamical properties of the
third-order correction models (\ref{GFM1}) in terms of K\"ahler or
para-K\"ahler geometry? 2) Except the Hermitian symmetric spaces
$G_{n,k}$ of compact type (i.e, {\it A III}), there are {\it C I},
{\it D III} and {\it BD I}-types of symmetric spaces (and also two
exceptional Hermitian symmetric spaces {\it E III} and {\it E VII})
(refer to \cite{He,FK}). Can we establish a similar result by using
the generalized bi-Schr\"odinger flow of maps from $\mathbb R$ to
these symmetric spaces? How about the other Hermitian symmetric
spaces of noncompact type and the other para-Hermitian symmetric
spaces? These questions deserve study in the future.

\section * {\bf Acknowledgements}
Part of this work was done when the first author visited SMC from
September to December, 2013. He thanks SMC for the invitation and
support. The authors are supported by the National Natural Science
Foundation of China (Q.D. by Grant No.10971030; Y.W. Grant
No.10990013),  RFDP (No. 20110071110002) and the innovation
foundation of Shanghai Educational Committee.


\begin{thebibliography}{50}
\setlength{\itemsep}{-3pt} \small
\bibitem{AKNS}
Ablowitz M.J., Kaup D.J., Newell A.C., Segur H.: The inverse
scattering transform-Fourier analysis for nonlinear problems. Stud.
in Appl. Math. 54, 249-315 (1974).

\bibitem{Ab0}
Ablowitz M.J., Clarkson P. A.: Nonlinear evolution equations and
inverse scattering. Cambridge Univ. Press, Cambridge (1992).

\bibitem{AMT}
Alekseevsky D., Medori C., Tomassina A.: Homogenous para-K\"ahler
Einstein manifolds. Russian Math. Surv. 64, 1-43 (2009).

\bibitem{AH}
Arms R.J., Hama F.R.: Localized-induction concept on a curved vortex
and motion of elliptic vortex ring. Phys. Fluids 8, 553-559 (1965).

\bibitem{AF}
Athorne C., Fordy A.P.: Generalized KdV and mKdV equations
associated with symmetric spaces. J. Phys. A 20,  1377-1386 (1987).

\bibitem{BJOS} Baldo S., Jerrard R.L., Orlandi G., Soner H.M.:
Vortex density models for superconductivity and superfluidity.
Comm. Math. Phys. 318, 131-171 (2013).


%\bibitem{AF}
%C. Athorne and A.P. Fordy, Generalized KdV and mKdV equations
%associated with symmetric soaces, J. Phys. A, {\bf 20} (1987)
%1377-1386.

\bibitem{BD}
Br\"ocker T., tom Dieck T.: Representations of compact Lie groups.
Springer-Verlag, New York, Berlin, Tokyo (1985).

\bibitem{BK}
Brower B.C., Kessler D.A., Koplik J.,Levine H.: Geometrical models
of interface evolution. Phys. Rev. A 29,  1335-1342 (1984).

%\bibitem{Ca}
%R. Calini, T. Ivey and G. Maris-Beffa, \emph{Remarks on KdV-type flows on
%star-shaped curves}, Phys. D {\bf 238} (2009) 788-797.


%\bibitem{CPe}
%S.S. Chern and C.K. Peng, \emph{Lie groups and KdV equations}, Manuscripta
%Math., {\bf 28} (1979), 207-217.

\bibitem{Chenb}
Chen B.B.: Schr\"odinger flows to symmetric spaces and the second
matrix-AKNS hierarchy. Comm. Theor. Phys. (Beijing) 45, 653-656
(2006).

\bibitem{CFG}
Cruceanu V., Fortuny P., Gadea P.M.: A survey on paracomplex
geometry. Rocky Mountain Journal of Mathematics 26, 83-115 (1996).

\bibitem{Da}
Da Rios L.S.: On the motion of an unbounded fluid with a vortex
filament of any shape, Rend. Circ. Mat. Palermo 22, 117-135 (1906).

\bibitem{D1}
Ding Q.: A note on the NLS and Schr\"odingr flow of maps. Phys.
Lett. A 248, 49-58 (1998).

\bibitem{ding2}
Ding Q.: Explicit blow-up solutions to the Schr\"odinger maps from
$\mathbb R^2$ to the hyperbolic 2-space $\mathbb H^2$. J. Math.
Phys. 50, no.10, 103507 (17pp) (2009).

\bibitem{DI}
Ding Q., Inoguchi J.: Schr\"odinger flows, binormal motion of curves
and the second AKNS hierarchies. Chaos, Solitons and Fractals 21,
669-677 (2004).

\bibitem{dinghe}
Ding Q., He Z.Z.: The Noncommutative KdV Equation and Its
Para-K\"ahler Structure. Sci. in China Math. to appear (2014).

\bibitem{DLW}
Ding Q, Liu X., Wang W.: The vortex fiament in the Minkowski
3-space. J. Phys. A 45, 455201 (14pp) (2012).

\bibitem{DWW}
Ding Q., Wang W., Wang Y.D.: A motion of spacelike curves in the
Minkowski 3-space and the KdV equation. Phys. Lett. A 374, 3201-3205
(2010).

\bibitem{DWW1}
Ding Q., Wang W., Wang Y.D.: The Fukumoto-Moffatt's model in the
vortex filament and generalized bi-Schr\"odinger maps. Phys. Lett. A
375, 1457-1460 (2011) .

\bibitem{DW1}
Ding Q., Wang Y.D.: Geometric KdV flows, motions of curves and the
third order system of the AKNS hierarchy. Inter. J. of Math. 22,
1013-1029 (2011).

\bibitem{ding1}
Ding Q., Zhu Z.: On the gaueg equivalent structure of the
Landua-Lishitz equation and its applications. J. Phys. Soc. of Japan
72, 49-53 (2003).

\bibitem{DW}
Ding W.Y., Wang Y.D.: Schr\"{o}dinger flows of maps into symplectic
manifolds. Sci. China  A 41, 746-755 (1998).

%\bibitem{DS}
%A. Doliwa and P.M. Santini; \emph{An elementary geometric characterization
%of the integrable motions of curves}, Phys. Lett. A, {\bf 185} (1994)
%373-384.

\bibitem{F}
Fukumoto Y.: Three-dimensional motion of a vortex filament and its
relation to the localized induction hierarchy. Eur. Phys. J. B 29,
167-171 (2002).

\bibitem{FM}
Fukumoto Y., Miyazaki T.: Three-dimensional distortions of a vortex
filament with axial velocity. J. Fluid. Mech. 22, 369-416 (1991).

\bibitem{FMo}
Fukumoto Y., Moffatt H.K.: Motion and expansion of a visous vortex
ring, I. A higer-order asymptotic formula for the velocity. J. Fluid
Mech. 417, 1-45 (2000).

\bibitem{FK}
Fordy A.P., Kulish P.: Nonlinear Schr\"odinger equations and simple
Lie algebras. Comm. Math. Phys. 80, 427-443 (1983).

%\bibitem{GP}
%R.E. Goldstein and D.M. Pertich; \emph{The Korteweg-der Vries hierarchy as
%dynamics of closed curves in the plane}, Phys. Rev. Lett., {\bf 67}
%(1991) 3203-329.

\bibitem{Hama}
Hama F.R.: Progressive deformation of a curved vortex filament by
its own induction. Phys. Fluids 5, 1156-1162 (1962).

\bibitem{Ha}
Hasimoto H.: A soliton on a vortex filament. J. Fluid. Mech. 51,
477-485 (1972).


\bibitem{He}
Helgason S.: Differential geometry, Lie group and symmetric spaces.
Academic Press, New York, San Francisco, London (1978).

\bibitem{Hel} Helmholtz H.: \"Uber Integrale der hydrodynamischen Gleichungen,
welche den Wirbelbewegungen entsprechen. J. Reine Angew. Math. 55,
25-55 (1858).

\bibitem{Hel*} Helmholtz H.: On integrals of the hydrodynamical equations, which
express vortex-motion (Translated by P. G. Tait). Phil. Mag., Ser.4,
33, 485-512  (1867).

\bibitem{Ke} Kelvin L.: The translatory velocity of a circular
vortex ring. Phil. Mag. 35, 511-512 (1867).

%\bibitem{K} B. Khesin; \emph{Dynamics of symplectic fluids and point vortices},
%Geom. Funct. Anal. 22 (2012), no. 5, 1444-1459.

%\bibitem{Hi} R. Hirota, Exact envelope-soliton solutions of a
%nonlinear wave equation, J. Math. Phys. \textbf{14} (1973), 805-809.

%\bibitem{KSVM}J. Klaers, J. Schmitt, F. Vewinger, M. Weitz;
%\emph{Bose-Einstein condensation of photons in an optical microcavity}, Nature 468(2010), 545¨C548.

%\bibitem{Ko}
%S. Kobayashi; \emph{Transformation groups in differential geometry},
%Berlin, Springer-Verlag, 1972.

%\bibitem{KL} J. Koplik and H. Levine; \emph{Vortex reconnection in
%superfluid helium}, Phys. Rev. Lett. 71(1993), 1375-1378.

\bibitem{Lan}
Langer J.S.: Instabilities and pattern formation in crystal growth.
Rev. Mod. Phys. 52, 1-28 (1980).

\bibitem{LP*} Langer J., Perline R.: Poisson geometry of the filament equation. J. Nonlinear
Sci. 1, 71-93 (1991).

\bibitem{LP}
Langer J., Perline R.: Geometric realizations of Fordy-Kulish
nonlinear Schr\"odinger sysyems. Pacific J. Math. 195, 157-178
(2000).

%\bibitem{MW} J. Marsden and A. Weinstein, \emph{Coadjoint orbits, vortices, and Clebsch variables
%for incompressible fluids}, Physica, 7D (1983), 305-323.


%\bibitem{MR} K. Moffatt, L. Ricca; \emph{Interpretation of invariants of the
%Betchov-Da Rios equation and of the Euler equations}, The Global
%Geometry of Turbulence, Plenum Press, New York 1991.

\bibitem{Lax}
Lax P.D.: Integrals of nonlinear equations of evolution and solitary
waves. Comm. Pure Appl. Math. 21, 467-490 (1968).


\bibitem{Liber}
Libermann P.: Sur les structures Presque paracomplex. C.R. Acad.
Sci. Paris 234, 2517-2519 (1952).

%\bibitem{LWY} T. Lin, J. Wei, J. Yang; \emph{Vortex rings for the Gross-Pitaevskii equation in $R^3$ },
%J. Math. Pures Appl. (9) 100 (2013), no. 1, 69-112.

%\bibitem{MCWD} K. W. Madison, F. Chevy, W. Wohlleben and J. Dalibard;
%\emph{Vortex Formation in a Stirred Bose-Einstein Condensate}, Phys. Rev. Lett. 84(2000), 806¨C809.
\bibitem{Ma}
Magri F.: A simple model of the integrable Hamiltonian equation. J.
Math. Phys.  19, 1156-1162 (1978).


\bibitem{MGK} Meleshko V.V., Gourjii A.A., Krasnopolskaya T.S.:
Vortex ring: history and state of the art. J. of Math. Sciences 187,
772-806 (2012).

\bibitem{OS}
Olver P.J., Sokolov V.V.: Integrable evolution equations on
associative algebra. Comm. Math. Phys.  193, 245-268  (1998).

\bibitem{ppetersen}
Petersen P.: Riemannian geometry . (2nd ed.) Springer Science
Business Media, LLC (2006).


\bibitem{PDL}
Porsezian K., Daniel M., Lakshmanan M. On the integrable models of
the higher one-dimensional classical continuum isotropic biquadratic
Heisenberg spin chain. J. Math. Phys. 33, 1807-1816 (1992).

\bibitem{Poh}
Pohlmeyer K.: Integrable Hamiltonian systems and interactions
through quadratic constraints. Comm. Math. Phys.  46, 207-221
(1976).

\bibitem{Ras}
Rashevskii P.K.: Scalar fields in a fibre space. Tr. Semin. Vekt.
Tenzor. Analiz. 6, 225-248 (1948).

%\bibitem{Ric}
%R.L. Ricca; \emph{Rediscovery of Da Rios equation}, Nature, {\bf 352}
%(1991), 561-562.

%\bibitem{ST}
%P.G. Saffman, G.I. Taylor, The penetration of a fluid into a porous
%medium or Hele-Shaw cell containing or more viscous liquid, Proc. R.
%Soc. London A {\bf 245} (1958) 312-329.

\bibitem{ST}
Saffman P.G., Taylor G.: The penetration of a fluid into a porous
medium or Hale-Shaw cell containing a more visous. Proc. R. Soc.
London A  245, 312-329 (1958).


%\bibitem{So}
%V.V. Sokolov and T. Wolf; \emph{Classiffication of integrable polynomial
%vector evolution equations}, J. Phys. A: Math. Gen., {\bf 34} (2001),
%11139-11148.

\bibitem{SW}
Sun X.W., Wang Y.D.: KdV geometric flows on K\"ahler manifolds.
Inter. J. Math. 22, 1439-1500 (2011).


\bibitem{Sym}
Sym A.: Soliton surfaces and their applications. Lecture Notes in
Physics  Vol.239, 145-231 (1985).


\bibitem{TeUh}
Terng C.L., Uhlenbeck K.: Schr\"odinger flows on Grassmannians. in
Integrable Systems, Geometry, and Topology, edited by Chuu-Lian
Terng, pp.235--256, AMS/IP Studies in Advanced Mathematics, Vol. 36,
(2006).

\bibitem{UIY}
Uby L., Isichenko L.B., Yankov V.V.: Vortex filament dynamics in
plasmas and superconductors. Physical Review E 52, 932-939 (1995).

\bibitem{ZT}
Zakharov V.E., Takhtajan L.A.: Equivalence of a nonlinear
Schr\"odinger equation and a Heisenberg ferromanet equation. Theor.
Math. Phys. 38, 17-23 (1979).
















\end{thebibliography}
\end{document}